\documentstyle{amsppt}

\magnification=1200
\NoBlackBoxes
\hsize=11.5cm
\vsize=18.5cm
\def\bt{\boxtimes}
\def\inv{^{-1}}
\def\dgla{differential graded Lie algebra\ }
\def\dg{differential graded\ }
\def\Aut{\text{\rm Aut}}
\def\CC{\Cal C}
\def\RR{\Cal R}
\def\M{\Cal M}
\def\E{\Cal E}
\def\SH{\Cal H}
\def\C{\Bbb C}
\def\H{\Bbb H}

\def\R{\Bbb R}

\def\lb{\lambda}

\def\sgn{\text{\rm sgn} }
\def\sym{\text{\rm sym} }
\def\Der{\text{\rm Der} }
\def\Der{\text{\rm Der}}
\def\Spec{\text{\rm Spec} }
\def\ad{\text{\rm ad}}
\def\bp{\bar{\partial}}
\def\ls{\vskip.25in}
\def\ss{\vskip.15in}

\def\O{\Cal O}
\def\UU{\Cal U}

\def\bt{\boxtimes}
\def\br{\text{br}}
\def\cob{^t\text{br}}
\def\Sym{\text{Sym}}
\def\id{\text{id}}
\def\rk{\text{rk}}
\def\Hi{\text{Hi}}
\def\hi{\text{hi}}
\def\g{\frak g}
\def\J{\frak J}
\def\gt{\tilde{\frak g}}
\def\gh{\hat{\frak g}}
\def\Et{\tilde{E}}
\def\Lt{\tilde{L}}
\def\gl{\frak{gl}}
\def\sll{\frak {sl}}
\def\h{\frak h}
\def\ab{\alpha\beta}
\def\m{\frak m}
\def\k{\frak k}
\def\endf{\frak { end}}
\def\hom{\frak {hom}}
\def\p{\partial}
\def\w{^\wedge}
\def\1/2{\frac{1}{2}}
\def\lmc{L_m(\C)}

\def\gsh{\frak g^\sharp}
\def\D{\frak D}
\def\U{\frak U}
\def\I{\Cal I}
\def\SS{\Cal S}
\def\Ann{\text{Ann}}
\def\nuline{\par\noindent}
\def\simto{\overset{\sim}\to{\to}}
\topmatter
\title {{Jacobi cohomology, local
geometry of moduli spaces, and Hitchin's connection}}\endtitle
\rightheadtext {Jacobi cohomology and Hitchin's connection}
\author Ziv Ran\endauthor
\address\nuline  Mathematics Department
University of California\nuline
Riverside CA 92521
USA\endaddress
\email   ziv\@math.ucr.edu\endemail
\date 8.8.1 \enddate
\toc

\head 0. Introduction \endhead
\head 1. Lie Atoms\endhead
\subhead 1.1 Basic notions\endsubhead
\subhead 1.2 Representations\endsubhead
\subhead 1.3 Universal enveloping atom\endsubhead
\head 2. Atomic deformation theory\endhead
\head 3. Universal deformations\endhead
\head 4. The Hitchin symbol\endhead
\subhead  4.1 The definition\endsubhead
\subhead 4.2 Cohomological formulae\endsubhead
\head 5. Moduli modules revisited\endhead
\head 6. Tangent algbera\endhead
\head 7. Differential operators\endhead
\head 8. Connection algebra\endhead
\head 9. Relative deformations over a global base\endhead
\head 10. The Atiyah class of a deformation\endhead
\head 11. Vector bundles on manifolds:
the action of base motions\endhead
\head 12. Vector bundles on Riemann surfaces:
refined action by base
motions and Hitchin's connection\endhead
\endtoc
\endtopmatter
\ls
\heading Introduction\endheading
The main purpose of this paper is to develop some
cohomological tools for the study of the local geometry of
moduli and parameter spaces in complex Algebraic Geometry.
The main ingredient will be
the language of Lie algebras, in particular differential graded
Lie algebras, their representations, and certain complexes
associated to these that we generally call
{\it{Jacobi}} complexes.
Why the presence of Lie algebras ? We understand since Felix
Klein that geometry, in one way or another, is conveniently
expressed in terms of symmetry groups, so it is reasonable to
expect a similar thing to be true of
deformations or variations
of a geometric object. Now a geometric
structure on a topological
space $X$ may be described
by gluing data on a collection of 'standard' or 'trivial'
pieces (e.g. polydiscs in the case of a manifold, or
or free modules in the case of a vector bundle),
and a deformation
of this structure may be obtained by varying
the gluing data.
Now, infinitesimal variations of gluing
data can be described
in terms of Lie algebras (e.g. of vector fields or
linear endomorphisms). Consequently,
infinitesimal deformations of geometric structures can
be systematically
expressed in terms of a sheaf of Lie algebras on $X$. Thus,
such sheaves will play a fundamental role in our work.\par
Actually, it often turns out to be convenient,
even necessary,
to work with a somewhat more general algebraic object than
Lie algebra, namely what we call a {\it {Lie atom}}.
Algebraically, a Lie atom is something like a quotient
of a Lie algebra by a subalgebra; to be precise, it consists
of a Lie algebra $\g$, a $\g-$module $\h$, together with a
module homomorphism $\g\to\h$. Geometrically, a Lie atom can
be used to control situations where a geometric object is
deformed while some aspect of the geometry 'stays the same'
(i.e. is deformed in a trivialized manner);
more particularly
the algebra $\g$ controls the deformation while the module
$\h$ controls the trivialization. Here we will present
a systematic development of some of the rudiments of
the deformation theory of Lie atoms,
which are closely analogous
to those of (differential graded) Lie algebras.\par
One of the main tools we develop here is a direct
cohomological construction, in terms of the moduli problem,
of vector fields and differential operators on moduli spaces,
together with their action on functions, as well as on
 'modular' modules, i.e. those associated
to the moduli problem,
 including formulae for composition and
Lie bracket (commutator); in
 particular, we obtain a canonical
formula for the Lie algebra
 of vector fields on a moduli space together with its natural
 representation on (formal) functions,
as well as extensions
to the case of differential operators acting on modular
vector bundles .\par
As an application of these methods we will
study the relation between
the geometry and deformations of a given complex manifold
$X$ and that of a moduli
space $\M_X$ of vector bundles on $X$. Since $\M_X$ is a functor
of $X$, it seems intuitively plausible that an automorphism
of $X$ should act on $\M_X$, and likewise for infinitesimal
automorphisms. This intuitive idea obviously
needs some precising,
because on the one hand the Lie algebra $T_X$ of holomprphic
vector fields on $X$ will typically admit no global sections,
and on the other hand as sheaves, $T_X$ and $T_{\M_X}$ live
on different spaces. In fact, we will show that there is a Lie
homomorphism $\Sigma_X$ from the \dgla
associated to $T_X$ to  that of
$T_{\M_X}$. This is useful because a Lie homomorphism
induces a map on the associated
deformation spaces, so $\Sigma_X$
can be used to relate deformations
of $X$ to those of $\M_X$.
\par
The latter result will be
further refined in case $X$ has dimension
1, i.e. is a compact Riemann surface, by showing that the map
$\Sigma_X$ factors through a Lie homomorphism to a certain Lie
atom associated to $\M_X$.
As an essentially immediate consequence
of this we will deduce the so-called Hitchin or Knizhnik-
Zamolodchikov flat connection over the moduli of curves.
This is a holomorphic connection on the projective bundle
associated to the vector bundle $\frak V$
with fibre $H^0(SU_X(r, L), G)$,
where $SU_X(r,L)$ is the moduli space of (S-equivalence
classes of) semistable bundles of rank $r$ and determinant $L$
on $X$, and $G$ is a line-bundle on $SU_X(r,L)$ (which
is necessarily, by results of
Drezet-Narasimhan [DrNa], a power of the
modular theta bundle, and a fractional power of the canonical
bundle). That the projectivization of $\frak V$ should
admit a flat connection was
conjectured by physicists based on ideas
from Conformal Field Theory,
and subsequently treated by a number
of mathematicians including
Beilinson-Kazhdan, Hitchin, Faltings,
Ueno and Witten (cf. [BeK] [BryM] [Hit] [Fa][Ram] [TsUY]
[vGdJ] [WADP] and
references therein). Our approach is quite close to Hitchin's
as regards the construction of the connection; the ideas here
go back to some degree to Welters [Wel]. However we are
able to extend the Welters-Hitchin construction, which
is essentially first-order deformation theory, to the Lie
theoretic context via what we call a {\it{connection algebra}},
which shows that the connection thus obtained is automatically
flat-- modulo showing that the
relevant maps are Lie homomorphisms.
We thus obtain a new and essentially 'algebraic'
proof of the flatness of the connection, replacing some arguments
by Hitchin [Hit] which appeal to
infinite-dimensional symplectic
geometry.\ls
The paper is organized as follows.
\S 1, 2 discuss basic definitions
and examples relating to
Lie atoms and their associated deformation
theory. In \S 3 we give a construction, under suitable
hypotheses, of the universal deformation associated to a
Lie atom, following closely the case of a Lie algebra.
In \S 4 we give a construction of the {\it{Hitchin symbol}}
attached to a family of curves, which is a crucial
ingredient in the contruction of the 'refined action'
by base vector fields on moduli spaces of vector bundles.
Whereas the usual construction of Hitchin symbols a in
[Hit, vGdJ] is based on Serre duality, hence is strictly
global on the curve, we realize the symbol as
the coboundary associated to a certain natural short
exact sequence, which later facilitates the proof of
some compatibilities with Lie brackets.\par
Next we revisit in \S 5 the construction of modular
modules, first given in [Ruvhs], and present it in
a new and more workable algebraic setting, based on
certain 'L complexes', which are 'adjoints' of the more
familiar modular Jacobi complexes. Based on this we give
in \S 6 the 'synthetic' construction of vector
fields and Lie brackets on moduli spaces, and in \S 7
the natural extension of these results to the case
of differential operators on modular modules.
\par
Next we present in \S 8 the notion of 'connection
algebra'. Given a Lie algebra $\g$ and a $\g-$module
$E$, the connection algebra $\k(\g,E)$ is a larger
\dgla having $\g$ as a quotient $E$ as a module.
It has the property that the cohomology of $E$
deforms in a trivialized way (i.e. carries a natural
flat connection) over the deformation space of $\k(\g,E)$.
This is a useful tool in the construction of Hitchin
connections.\par
In \S 9 we discuss the extension of the foregoing results
to the case  of {\it{relative}} deformations. Whereas
an ordinary deformation is considered parametrized by
a thickened {\it{point}} (i.e. an artin local algebra),
a relative deformation is likewise parametrized by a
thickened {\it{space}} (i.e. a coherent algebra over a
ringed space). In \S 10 we discuss an analogue, in the
setting of deformation theory, of the notion of {\it{Atiyah
class}} or {\it{Atiyah extension}}. We show that
the Atiyah extension is an extension of Lie algebras
admitting a natural representation.\par
The final two sections focus on applications of the foregoing
techniques to moduli spaces of vector bundles on a
manifold and their deformation spaces. In \S 11
we construct the map $\Sigma_X$ mentioned above as
a homomorphism of \dg Lie algebras, which gives us
a precise handle on the relation between deformations
of $X$ and those of $\M_X$. Then in \S 12 we construct,
based in the Hitchin symbol as presented in \S 4,
a lifting of $\Sigma_X$ (as Lie homomorphism) to
a certain Lie atom associated to $\M_X$, which is
related to a suitable connection algebra.
This yields the Hitchin connection and its flatness
essentially for free.\ls
Some of the constructions and techniques in this paper
are presented in greater generality than is required just
for the Hitchin connection. Hopefully they may find other
applications to the geometry of moduli spaces.\ls
\remark{Acknowledgment} Some of the work on this paper was done
while the author was visiting the Mathematics Department at
Roma Tre University. He is grateful to the department,
and especially to Angelo Lopez, Edoardo Sernesi and Sandro Verra,
for their hospitality and for providing a congenial
and stimulating working environment.\endremark
\vfill\eject
\heading 1.  Lie atoms\endheading

Our purpose here is to
define and begin to study a notion
which we call {\it{Lie atom}}
and which generalizes that of the quotient of a Lie
algebra by a subalgebra (more precisely, a pair of Lie
algebras viewed in the derived category).
Our point of view is that a Lie atom,
though not actually a Lie algebra, possesses
some of the formal properties of Lie algebras. In particular,
we shall see later that there is a deformation theory for Lie atoms,
which generalizes the case of  Lie algebras  and which in addition
allows us to treat some classical, and disparate,
 deformation problems such as, on the one hand,
the Hilbert scheme, and on the other hand heat-equation deformations,
introduced in the first-order case by Welters [We].

 \subheading{1.1 Basic notions}
\proclaim{Definition 1.1.1} By a {\rm{ Lie atom}}
(for 'algebra to module') we shall mean the data
$\gsh$ consisting of
\item{(i)} a Lie algebra $\g$;
\item{(ii)} a $\g$-module
$\h$;
\item{(iii)} an injective $\g$-module homomorphism
$$i: \g\to\h,$$
where $\g$ is viewed as a $\g$-module via the adjoint action.
\endproclaim
The assumption of injectivity is not really essential but is
convenient and is satisfied in applications.
Hypothesis (iii) means explicitly that,
writing $<\ ,\ >$ for the $\g$-action on $\h$, we have
$$i([a,b])=<a,i(b)>=-<b,i(a)>.\tag 1.1$$
Note that any Lie algebra $\g$ determines a 'Lie atom', minus
the injectivity hypothesis, by taking
$\h=0,$  and the concept of Lie atom is essentially
a generaization of that of Lie algebra.
Note also that there is an obvious notion of morphism
of Lie atoms, hence also of isomorphism and quasi-isomorphism
(composition of morphisms inducing isomorphism on cohomology and
inverses of such).
Of course one can also talk about sheaves of Lie atoms,
\dg\ Lie atoms, etc. We shall generally
consider two atoms to be equivalent
if they are quasi-isomorphic.
\remark{ Examples 1.1.2}\par{\bf A.}
If $j:E_1\to E_2$ is any linear map of
vector spaces, let $\g=\g(j)$ be the
{\it{interwining algebra}} of $j$,
i.e. the Lie subalgebra
$$\g\subseteq\gl(E_1)\oplus\gl(E_2)$$ given by
$$\g =\{ (a_1,a_2)|j\circ a_1=a_2\circ j\} .$$
Thus $\g$ is the 'largest' algebra
acting on $E_1$ and $E_2$
so that $j$ is a $\g-$homomorphism.
When $j$ is injective, define
$$\gl (E_1<E_2):=(\g, \gl (E_2), i),$$
with $\ \ i(a_1,a_2)=a_2.\ \ $ When $j$ is surjective, define
$$\gl (E_1>E_2):=(\g, \gl (E_1), i),$$
with $\ \ i(a_1,a_2)=a_1 \ \ $
These are Lie atoms. Again, the
definitions could be made without
assuming $j$ injective or surjective, but we have no
interesting examples. The two notions
are obviously dual to each other, but since we do not assume
$E_1, E_2$ are finite-dimensional, dualising is not necessarily
convenient.\par
{\bf B.} If $i:\g_1\to\g_2$ is an injective
  homomorphism of Lie algebras then
$$\gsh :=(\g_1,\g_2,i)$$
is a Lie atom. More generally,
if $\h$ is any $\g_1$ submodule of
$\g_2$ containing $i(\g_1)$, then
$$\gsh :=(\g_1,\h,i)$$
is a Lie atom.\par
 {\bf C.} Let $E$ be an invertible sheaf on a ringed space $X$
(such as a real or complex manifold), and let
$\D^i(E)$ be the sheaf of $i-$th order differential endomorphisms
of $E$ and set
$$\D^{\infty}(E)=\bigcup\limits_{i=0}^{\infty}\D^i(E).$$
Then $\g=\D^1(E)$ is a Lie algebra sheaf and $\h=\D^2(E)$
is a $\g-$module, giving rise to a Lie atom $\gsh$ which will be called
the {\it{Heat atom}} of $E$ and denote by $\D^{1/2}(E)$.
Note that if $X$ is a
manifold then $\gsh$ is quasi-isomorphic as a complex to
$\Sym^2(T_X)$. \par
{\bf D.} Let $$Y\subset X$$
be an embedding of manifolds (real or complex). Let
$T_{X/Y}$ be the sheaf of vector fields on $X$ tangent
to $Y$ along $Y$. Then $T_{X/Y}$ is a sheaf of Lie algebras
contained in its module $T_X$, giving rise to a Lie atom
$$N_{Y/X}= (T_{X/Y}\subset T_X),$$
which we call the {\it{normal atom}} to $Y$ in $X$. Notice that
$T_{X/Y}\to T_X$ is locally an isomorphism off $Y$, so replacing
$T_{X/Y}$ and $T_X$ by their sheaf-theoretic restrictions on $Y$
yields a Lie atom that is quasi-isomorphic to, and identifiable with
$N_{Y/X}$.\par
{\bf E} In the situation of the previous example, let $\I_Y$ denote
the ideal sheaf of $Y$. Then $\I_Y.T_X$ is also a Lie subalgebra
of $T_X$
giving rise to a Lie atom $$T_X\otimes\O_Y:=(\I_Y.T_X\subset T_X).$$
Note that via the embedding of $Y$ in $Y\times X$ as the graph
of the inclusion $Y\subset X$, $T_X\otimes\O_Y$ is quasi isomorphic as
Lie atom to $N_{Y/Y\times X}$, so this example is essentially a
special case of Example D.
\endremark\ls
\subheading{1.2 Representations}
Now given a Lie atom $\gsh=(\g,\h,i)$, by a {\it{left $\gsh-$module}}
or {\it{left $\gsh-$representation}}  we shall mean
the data of a pair $(E_1,E_2)$ of $\g-$modules with an injective $\g-$
homomorphism $j:E_1\to E_2$, together with an 'action rule'
$$<\ \ >:\h\times E_2\to E_2,$$
satisfying the compatibility condition (in which we have written
$<\ >$ for all the various action rules):
$$<<a,v>,x>=<a,<v,x>>-<v,<a,x>>,\tag 1.2$$
$$\forall a\in\g , v\in\h ,x\in E_1.$$
In other words, a left $\gsh-$ module is just a homomorphism of Lie atoms

$$\gsh\to\gl (E_1<E_2).$$
The notion of right $\gsh-$module is defined similarly.

\remark{Examples, bis}  Refer to the previous
examples.\item {\bf A.} These are the tautological examples:
$\gl(E_1<E_2)$ and $\gl(E_1>E_2)$ with $(E_1,E_2)$
as left (resp. right) module in the two cases $j$ injective
(resp. surjective).\item
{\bf B.}For a Lie atom $\gsh=(\g_1,\g_2,i)$, $\gsh$ itself is a
left  $\gsh$-module,
called the {\it{adjoint representation}}
while $(\gsh)^*=(\g_2^*, \g_1^*, i^*), *=$ dual
vector space, is a right $\gsh-$module called
the {\it{coadjoint representation}}.\item
{\bf C.} In this case $(E,E)$ is a left and right
$\gsh-$module, called a {\it{Heat module}}.\item
{\bf D.} Here the basic left module is $(\I_Y,\O_X)\overset{qis}\to\sim
\O_Y.$ Of course we may replace $\I_Y$ and $\O_X$ by
their topological restrictions of $Y$. The basic right module
of interest is $(\O_X,\O_Y)$.\item
{\bf E.} In this case the modules we are interested in are
$$(\I_{Y,Y\times X},\O_{Y\times X}),\ \ (\O_{Y\times X}
\to \O_Y).$$
\endremark
\remark{Remark} 'Theoretically', only the action of $\h$
going from $E_1$ to $E_2$ 'should' be necessary
for a module. However the action
on all of $E_2$ is needed in proofs and satisfied in the examples
we have in mind, so we included it. The fact that we require an
extension to $E_2$ rather than the dual notion of
a 'lifting' to $E_1$ has to do with
the fact that in our examples of interest the maps $i$ and $j$
are injective.
\endremark
\subheading{1.3 Universal enveloping atom}
Observe that for any Lie atom $(\g,\h,i)$ there is a smallest
Lie algebra $\h^+$ with a $\g-$ map
$\h\to\h^+$  such that the given action
of $\g$ on $\h$ extends via $i$ to a 'subalgabra'
action of $\g$ on $\h^+$, i.e. so that
$$<a,v>=[i(a), v],\ \forall a\in\g, v\in\h^+,\tag 1.3$$
namely $\h^+$ is simply the quotient of the
free Lie algebra on $\h$ by the ideal generated by elements of
the form
$$[i(a), v]- <a,v>,\ \ a\in\g, v\in\h$$
(note that the action of $\g$ on $\h$ extends to an action
on $\h^+$ by the 'derivation rule'). In view of the basic
identity (1.1) it follows that the map $\g\to\h^+$ induced by $i$
is a Lie homomorphism. This shows in particular that,
modulo replacing $\g$ and $\h$ by their images in
$\h^+$, any Lie
atom $\gsh$ is essentially of the type of Example B above
(though the map $\h\to\h^+$ is not necessarily injective).\par
We observe next that there is an natural notion of 'universal
enveloping atom' associated to a Lie atom $\gsh=(\g,\h,i)$.
Indeed let $\U(\g,\h)$ be the quotient of the
$\U(\g)$-bimodule $\U(\g)\otimes\h\otimes\U(\g)$ by the sub-bimodule
generated by elements of the form
$$a\otimes v-v\otimes a-<a,v>, i(a)\otimes b-a\otimes i(b),$$
$$\forall a,b\in\g, v\in\h.$$
\noindent{\bf{Sorites }}\par
{1.}\  $\U(\g,\h)$ is a $\U(\g)-$bimodule .\par
2.\ The map $i$ extends to a bimodule homomorphism
$$i:\U(\g)\to\U(\g,\h).$$\par
{3.}\   $\U(\g,\h)$ is universal with
respect to these properties.\par
4.\ $\U(\g,\h)$ is generated by $\h$ as either right or left
 $\U(\g)-$module. Moreover the image of $\U(\g,\h)$
 in $\U(\h^+)$ is precisely
 the (left, right or bi-)$\U(\g)$-submodule of $\U(\h^+)$
 generated by $\h$. \par
\ls
Thus $$\U(\gsh):=(\U(\g),\U(\g,\h),i)$$ forms an 'associative atom' which we
call the {\it{universal enveloping atom}} associated to $\gsh$.
\remark{Examples, ter }\item
{\bf A.} It is elementary that the universal enveloping algebra
of the interwining Lie algebra $\g$ is simply the interwining associative
algebra
$$\U(\g)=\{(a_1,a_2)|j\circ a_1=a_2\circ j\},$$
and so the universal enveloping atom of $\gl (E_1<E_2)$
(resp. $\gl (E_1>E_2)$ is just
$$(\U(\g), \endf (E_2), i)\ \ resp.\ \ (\U(\g), \endf (E_1), i) .$$\item
{\bf B. } In this case it is clear that $\U(\g_1,\h)$
is just the sub $\U(\g_1)-$bimodule generated by $\h$.
\endremark
\vfill\eject
\heading 2. Atomic deformation theory\endheading
Our purpose here is to define and study deformations with
respect to a Lie atom $\gsh=(\g,\h,i)$.
Roughly speaking a $\gsh-$deformation consists of a  $\g-$
deformation $\phi$, plus a 'trivialization of $\phi$ when
viewed as $\h-$ deformation;'  as we shall see in the course
of making the latter precise, it only involves the structure of $\h$
as $\g-$module, not as Lie algebra, and this is our main
motivation for introducing the notion of Lie atom.\par
We recall first the notion of
$\g-$deformation. Let $\g$ be a sheaf of Lie algebras over
a Hausdorff topological space $X$, let $E$ be a $\g-$module
and $S$ a finite-dimensional $\C-$algebra with maximal ideal
$\m$. Note that there is a sheaf of groups
$G_S$ given by
$$G_S=\exp (\g\otimes\m)$$
with multiplication given by the Campbell-Hausdorff formula,
where $\exp$, as a map to $\U(\g\otimes\m)$, is injective because
the formal $\log$ series gives an inverse.
Though not essential for our purposes, it may be noted
that $G_S$ coincides with the (multiplicative)
subgroup sheaf of sections congruent to 1 modulo the ideal
$\U^+ (\g\otimes \m)$ generated
by $\g\otimes\m$ in the universal
enveloping algebra
$$\U(\g\otimes \m):=\U_S(\g\otimes \m).$$ This is easy to prove by induction on
the exponent of $S$: note that if $I<S$ is an ideal with
$\m .I=0$ then $\g\otimes I\subseteq \g\otimes\m$ is a central ideal
yielding a central subgroup $G_I=1+\g\otimes I\subseteq G_S$ and a central ideal
$\g\otimes I\subseteq\U(\g\otimes\m)$, hence a central
subgroup $1+\g\otimes I$ in the multiplicative group of
$\U(\g\otimes\m).$ \par
 A $\g-$deformation of $E$ over $S$ is a sheaf $E^{\phi}$
of $S$-modules, together with a maximal atlas of trivialisations
$$\Phi_{\alpha}:E^{\phi}|_{U_{\alpha}}\overset{\sim}\to{\to}
 E|_{U_{\alpha}}\otimes S,$$
such that the transition maps
$$\Psi_{\alpha\beta}:=\Phi_{\beta}\circ\Phi_\alpha^{-1}\in
G_S(U_{\alpha}\cap U_{\beta}).$$ We view a $\g-$deformation
(not specifying any $E$) as being
given essentially by the class of $(\Psi_{\alpha\beta})$ in the
nonabelian \v{C}ech cohomology set $\check{H}^1(X,G_S)$ and in
particular a $\g-$deformation determines simultaneously
$\g-$deformations of all $\g-$modules $E$, and is in turn
determined by the corresponding $\g-$deformation of
any {\it{faithful}} $\g-$module $E$ . We may call $E^\phi$
a model of $\phi$ or $(\Psi_{\alpha\beta})$.\par
 Now let
$\gsh=(\g,\h,i)$ be a sheaf of Lie atoms on $X$,
and let $E^\sharp =(E_1,E_2,j)$
be a sheaf of left $\gsh-$modules.
Note that an element $v\in\h\otimes\m$ determines
a map
$$A_v:E_1\to E_2,$$
$$A_v(x)=j(x)+<v,x>.\tag 2.1$$
Locally, an $S-$linear map
$$A:E_1\otimes S\to E_2\otimes S$$
is said to be a {\it{left $\h-$map}} if it is of the form
$$A=\exp(u)\circ A_v, u\in\g\otimes \m, v\in \h\otimes\m,$$
and similarly for right modules and right $\h-$maps.
We note that the set of left (resp. right) $\h-$maps is invariant under
the left (resp. right) action of $G_S$ on $\hom (E_1\otimes S,E_2\otimes S)$.
We consider the data of an $\h-$map (left or right) to include
the element $v\in\h\otimes\m$, and two such maps are considered
equivalent if they belong to the same $G_S-$orbit. Thus a left
$\h-$ map is really a $G_S-$orbit in $G_S.(\h\otimes\m).$
\par
The notion of $\h-$map
globalizes as follows. Given a  
 $\g-$deformation $\phi = (\Psi_{\alpha\beta})$,
a (global) {\it{left $\h-$map}} (with respect to $\phi$) is a map
$$A:E_1\otimes S\to E_2^{\phi}$$
such that for any atlas
$\Phi_\alpha$ for $E_2^{\phi}$
over an open covering $U_\alpha$,
$\Phi_\alpha\circ A$ is given over $U_\alpha$ by a left $\h-$map.
Note that this condition is independent of the choice of atlas,
and is moreover equivalent to the {\it existence} of some
atlas for which the $\Phi_\alpha\circ A$ are given by
$$x\mapsto j(x)+<v_\alpha,x>,\ \  v_\alpha\in\h (U_\alpha)\otimes \m.
\tag 2.2$$
We call such an atlas a {\it{good atlas}} for $A$.
The notion of global right $\h-$map
$$B:F_1^\phi\to F_2\otimes S$$
for a right $\gsh-$module $(F_1, F_2, k)$
is defined similarly, as are those of global left and right
$\h-$maps without specifying a module. A pair $(A,B)$ consisting
of a left and right $\h-$map is said to be a {\it{dual pair}}
if there exists a common good atlas with respect to which
$A$ has the form (2.2) while $B$ has the form
$$x\mapsto j(x)-<v_\alpha,x>$$
with the {\it{same}} $v_\alpha$.
\proclaim{Definition} In the above situation, a left $\gsh-$deformation of
$E^\sharp$ over $S$ consists of a $\g-$deformation $\phi$
together with an  $\h-$map  from the trivial
deformation to $\phi$:
$$A:E_1\otimes S\to E_2^{\phi}.$$
Similarly for right $\gsh-$deformation.
A $\gsh-$deformation consists of a $\g-$deformation $\phi$ together with
a dual pair $(A,B)$ of $\h-$maps with respect to $\phi$.

\endproclaim
In terms of \v{C}ech data $(\Psi_{\alpha\beta}=\exp (u_{\alpha\beta}), v_\alpha)$
for a good atlas as above, the condition that the $j+v_\alpha$
should glue together to a globally defined map left $\h-$ map $A$ is
$$\Psi_{\alpha\beta}\circ (i+v_\alpha )=i+v_\beta,\tag 2.3$$
which is equivalent to the following equation in $\U(\g\otimes\m,\h\otimes\m)$,
in which we set
$$D(x)=\frac{\exp (x)-1}{x}=
\sum\limits_{k=1}^\infty \frac{x^k}{(k+1)!}\ \  :\tag 2.4$$
$$D(u_{\alpha\beta})i(u_{\alpha\beta})+\exp (u_{\alpha\beta}).v_\alpha
=v_\beta . \tag 2.5$$ The condition for a right $\h -$ map $B$ is
$$i(u_{\alpha\beta})D(u_{\alpha\beta})+v_\alpha\exp (u_{\alpha\beta})
=v_\beta . \tag 2.6$$
\remark{Examples} We continue with the examples of \S 1.\item
{\bf C} When $\gsh = (\D^1(E),\D^2(E))$ is the heat algebra of
the locally free sheaf $E$,
$\gsh-$deformations of $E^\sharp=(E,E)$ are called {\it{heat deformations}}.
Recall that a $\D^1(E)-$ deformation consists of a deformation
$\O^\phi$ of the structure sheaf of $X$, together with an invertible
$\O^\phi-$ module $E^\phi$ that is a deformation of $E$. Lifting
this to a $\gsh-$ deformation amounts to constructing $S-$linear,
globally defined maps ( {\it{heat operators}})
$$A:E\otimes S\to E^\phi,$$
$$B:E^\phi\to E\otimes S$$
that are locally (with respect to an atlas and a trivialisation of $E$) of
the form
$$(f_i)\mapsto (f_i\pm\sum a_{j,k}\partial f_j/\partial x_k\pm\sum a_{j,k,m}
\partial^2f_j/\partial x_k\partial x_m),\tag 2.7$$
$$a_{j,k}, a_{j,k,m}\in\m\otimes\O_X.$$
Notice that the heat operator $A$ yields a well-defined lifting of sections
(as well as cohomology classes, etc.) of $E$ defined in any open set $U$
of $X$ to sections of $E^\phi$ in $U$. In particular,
suppose that $X$ is a compact complex manifold and that
$$H^i(X,E)=0,\ \ i>0.$$
It follows, as is well known, that $H^0(E^\phi)$ is a free
$S$-module, hence
$$H^0(A):H^0(X,E)\otimes S\to H^0(X,E^\phi)$$
is an isomorphism.
Thus for any heat deformation the module of sections
$H^0(E^\phi)$ is {\it{canonically trivialised}}.
Put another way,  $H^0(E^\phi)$ is endowed with
a {\it{canonical flat connection}}
$$\nabla^\phi:H^0(E^\phi)\to H^0(E^\phi)\otimes\Omega_S\tag 2.8$$
determined by the requirement that
$$\nabla^\phi\circ H^0(A) =0,\tag 2.9$$
i.e. that the image of heat operator consist of
flat sections. \item
{\bf D.} When $\gsh = N_{Y/X}$, a left $\gsh-$ deformation
of $(\I_Y,\O_X)$ consists of
a $T_{X/Y}-$deformation, i.e. a deformation
$(X^\phi,Y^\phi)$ of $(X,Y)$ in
the usual sense, together with  $T_X-$maps
$$A:\I_Y^\phi\O_X^\phi\to \O_X\otimes S,$$
$$B:\O_X\otimes S\to\O_X^\phi\to\O_Y^\phi.$$
which yield trivialisations   of the deformation $X^\phi$.
Thus left $\gsh-$deformations yield deformations of $Y$
in a fixed $X$, and similarly for right deformations. Conversely,
given a deformation of $Y$ in a fixed $X$, let $(x_\alpha^k)$
be local equations for $Y$ in $X$, part of a local
coordinate system. Then it is easy to see that
we can write equations for the deformation of $Y$ in the form
$$x_\alpha^k+v_\alpha (x_\alpha^k),\ \ v_\alpha\in T_X\otimes\m$$
($v_\alpha$ independent of $k$), so this comes from a left
and a right $\gsh-$ deformation of the form $((\Psi_{\alpha\beta}),(v_\alpha))$
where
$$\Psi_{\alpha\beta}=(1+v_\alpha)(1+v_\beta)^{-1}\in \U_S(T_{X/Y}\otimes\m)
(U_{\alpha}\cap U_\beta).$$
Thus the three notions of left and right $N_{Y/X}-$deformations and
deformations of $Y$ in a fixed $X$ all coincide.
\item
{\bf E.} In this case we see similarly that $T_X\otimes\O_Y-$ deformations
of $(\I_Y\oplus\O_Y,\O_X)$ consist of a deformation of the pair $(X,Y)$,
together with trivialisations of the corresponding deformations of $X$
and $Y$ separately, i.e. these are just deformations of the embedding
$Y\hookrightarrow X$, fixing both $X$ and $Y$.

\endremark
\vfill\eject
\heading 3. Universal deformations\endheading
Our purpose here is to construct the universal deformation
associated to a sheaf $\gsh$ of Lie atoms which is simultaneously
the universal $\gsh-$deformation of any $\gsh-$module $E^\sharp$.
We thus extend the main result of [Rcid] to the cases of atoms.
We refer to [Rcid] and [Ruvhs] for details on any items not explained here.
\par
We shall assume throughout, without explicit mention, that all
sheaves of Lie algebras and modules considered are {\it admissible}
in the sense of [Rcid].
In addition, unless otherwise stated we shall assume their
cohomology is finite-dimensional.
We begin by reviewing the main construction
of [Rcid] and restating its main theorem in a slightly stronger form.
The {\it{Jacobi complex}} $J_m(\g)$ is a complex in degrees
$[-m,-1]$ defined on $X<m>$,
the space parametrizing nonempty subsets of $X$ of cardinality at
most $m$ This complex has the form
$$\lb^m(\g)\to ...\to\lb^2(\g)\to \g$$
where $\lb^k(\g)$ is the exterior alternating tensor power
and the coboundary maps
$$d_k:\lb^k(\g)\to\lb^{k-1}(\g)$$
are given by
$$d_k(a_1\times ...\times a_k)=\frac{1}{2}\frac{1}{k!}\sum\limits_{\pi\in S_k}
[\sgn (\pi )[a_{\pi (1)},a_{\pi (2)}]\times a_{\pi (3)}...\times a_{\pi (k)}.
\tag 3.1$$
(This differs from the formula in [Rcid] by a factor of  1/2, which obviously
makes no essential difference but is convenient.) We showed in [Rcid ] that
$$R_m(\g)=\C\oplus\H^0(J_m(\g))^*$$
is a $\C-$algebra (finite-dimensional
by the admissibility hypothesis)
and we constructed a certain 'tautological'
$\g-$deformation $u_m$ over it.
To any $\g-$deformation $\phi$ over an algebra $(S,\m)$ of exponent $m$
 we associated a canonical Kodaira-Spencer
homomorphism
$$\alpha =\alpha (\phi ) :R_m(\g )\to S \ .$$
Although in [Rcid] we made the hypothesis that $H^0(\g)=0,$ this is in fact
not needed for the foregoing statements, and is only used in the proof that
$u_m$ is an $m-$universal deformation.\par
Now the hypothesis $H^0(\g)=0$ can be relaxed somewhat. Let us say that
$\g$ has {\it{central sections}} if for each open set
$U\subset X$, the image
of the restriction map
$$H^0(\g)\to \g(U)$$
is contained in the center of $\g(U)$.
Equivalently, in terms of a soft dgla resolution
$$\g\to\g^.,$$
the condition is that $H^0(\g)$ be contained in the center of
$\Gamma (\g^.)$, i.e. the bracket
$$H^0(\g )\times \Gamma (\g^i)\to\Gamma (\g^i)$$
should vanish.
\proclaim{Theorem 3.1} Let $\g$ be an admissible dgla
and suppose that
$\g$ has central sections. Then for any $\g-$deformation $\phi$
there exist isomorphisms
$$\phi\simto\alpha (\phi )^*(u_m)=u_m\otimes_{R_m(\g)}S,\tag 3.2$$
any two of which differ by an element of
$${\text{ Aut}} (\phi )=H^0(\exp (g^\phi\otimes\m)).$$
In particular, if $H^0(\g)=0$ then the isomorphism is unique.
Consequently, for any admissible pair $(\g,E)$ there are
isomorphisms
$$E^\phi\simto\alpha (\phi )^*(E^{u_m})$$
any two of which differ by an element of
${\text{ Aut}} (\phi ).$

\endproclaim
\demo{proof} This is a matter of adapting the argument in the proof
of Theorem 0.1, Step 4, pp.63-64 in [Rcid], and we will just
indicate the changes. We work in $\H^0(J_m(\g)) \otimes\m$
rather than $\H^0(J_m(\g), \m^.)$, which may not inject to it.
Then, with the notation of {\it{loc. cit.}} we may write
$$u_1=\sum v_i\otimes\phi_i\in \Gamma (\g^0)\otimes
\Gamma (\g^1)\otimes\m .$$
The argument there given shows that
$$u_1=u_0\times \phi + w\times\phi$$
where $w\times\phi\in H^0(\g)\otimes\Gamma (\g^1)\otimes \m.$
Now- and this is the point- since
$$[w,\phi]=0$$
thanks to our assumption of central sections, this is sufficient
to show that $\tilde{\phi}-\phi$ is the total coboundary of
$u_0\times \phi$, as required.
\par
 This shows the existence of the isomorphism
as in (3.2). Given this, the fact that two such isomorphisms
differ by an element of ${\text{ Aut}} (\phi )$ is obvious.
To identify the latter group it suffices to identify
its Lie algebra,
which is given by the set of $\g-$endomorphisms
$$\ad (v)\in\g^0\otimes\m$$ of the
resolution
$$(\g^.\otimes S, \bp +\ad (\phi )).$$
It is elementary to check that the condition on $v$ is precisely
$$\bp (v) +\ad (\phi )(v) =0,$$
i.e. $v\in\g^\phi\otimes\m$.
Finally since $\g^\phi$ is isomorphic as a sheaf
(not $\g-$isomorphic) to the trivial deformation
$\g\otimes S$, we have
$$H^0(\g^\phi)=H^0(\g\otimes\m)=H^0(\g)\otimes\m, $$
hence ${\text{ Aut}} (\phi )=(1)$ if $H^0(\g)=0.$
\qed

\enddemo
\remark{Remark} Without the hypothesis of central sections
it is still possible to 'classify' $\g-$deformations over
$(S,\m)$ in terms of $\H^0(J_m(\g), \m^.)$ but it is not
immediately clear how this is related to semiuniversal deformations.
We hope to return to this elsewhere.
\endremark
We now extend these results from Lie algebras to Lie atoms.
First recall the modular Jacobi complex $J_m(\g,\h)$ associated
to a $\g-$module $\h$. This is a complex in degrees $[-m,0]$
defined on a space $X<m,1>$ parametrizing pointed subsets of $X$
of cardinality at most $m+1.$ It has the form
$$\lb^m(\g)\bt\h\to ...\to\g\bt\h\to \h$$
with differentials
$$\partial_k(a_1\times ...\times a_k\times v)=d_k(a_1\times ...\times a_k)
\times v +\frac{1}{2k}\sum\limits_{j=i}^k (-1)^ja_1\times ...\hat{a_j}\times
...\times a_k\times (a_j(v)),$$
where $d_k$ is the differential in $J_m(\g)$ (see [Ruvhs]).\par
Now let $\gsh=(\g,\h,i)$ be a Lie atom. Then the $\g-$homomorphism
$i$ gives rise to a map of complexes
$$i_m:J_m(\g)\to \pi_{m-1,1\ *}J_{m-1}(\g,\h),$$
where $ \pi_{m-1,1}:X<m-1,1>\to X<m>$ is the natural map.
We denote by $J_m(\gsh)$ the mapping cone of $i_m$ and call it
the Jacobi complex of the atom $\gsh$.
We note the natural map
$$\sigma_m :J_m(\gsh)\to (J_m/J_1)(\gsh)\to\Sym^2(J_{m-1}(\gsh))$$
obtained by assembling together various 'exterior comultiplication'
maps
$$\lb^j(\g)\to \lb^r(\g)\bt\lb^{j-r}(\g),$$
$$\lb^j(\g)\bt\h\to\lb^r(\g)\bt (\lb^{j-r}(\g)\bt\h).$$
This gives rise a (commutative, associative) OS
(i.e. comultiplicative) structure on $J_m(\gsh)$, which
induces one on $H^0(J_m(\gsh))$, whence a local finite-dimensional
$\C-$algebra structure on
$$R_m(\gsh):= \C\oplus H^0(J_m(\gsh))^*$$
as well as a local homomorphism
$$R_m(\g)\to R_m(\gsh),$$
dual to the 'edge homomorphism'
$J_m(\gsh)\to J_m(\g)$.\par
\remark{Remark}If $\h$ happens to be a Lie algebra
and $i$ a Lie homomorphism
we may think of $R_m(\g/\h)$ as a formal functorial
substitute for
the fibre of the
induced homomorphism $R_m(\h)\to R_m(\g)$, i.e.
$R_m(\g)/\m_{R_m(\h)}. R_m(\g)$. It is important to
note that this fibre involves only the $\g-$module
structure on $\h$ and not the full algebra structure.
\endremark
Now we may associate to a $\gsh-$deformation a $J_m(\gsh)-$cocycle
as follows. By definition, we have a pair of pairs of
$\g-$deformations with a dual pair of $\h-$ maps
$$A: E_1\otimes S\to E_2^\phi,$$
$$B: E_3^\phi\to E_4\otimes S.$$
As usual we represent $E_i^\phi$ by a resolution of the from
$(E_i^.\otimes S, \bp +\phi), i=2,3$
and $E_i\otimes S$ by $(E_i^.\otimes S, \bp).$
Then the maps $A,  B$ can be represented
simultaneously in the form $j\pm v, v\in \h^1\otimes\m_S,$
and we get a pair of commutative diagrams
$$\matrix E_1^i\otimes S&\overset\bp\to\to&E_1^{i+1}\otimes S\\
j+v\downarrow&&\downarrow j+v\\
E_2\otimes S&\overset{\bp+\phi}\to\to&E_2^{i+1}\otimes S\endmatrix\tag 3.3$$
$$\matrix E_3^i\otimes S&\overset{\bp+\phi}\to\to&E_2^{i+1}\otimes S\\
j-v\downarrow&&\downarrow j-v\\
E_4\otimes S&\overset{\bp}\to\to&E_4^{i+1}\otimes S\endmatrix\tag 3.4$$
whose commutativity amounts to a pair of identities in $\U(\g,\h)$:
$$\phi .v =-\bp (v) -i(\phi),\ \ v.\phi=\bp(v)+i(\phi)\tag 3.5$$
which imply
$$\bp(v)+i(\phi)=-\1/2 <\phi,v>.\tag 3.6$$
Now the identity (3.6) together with the usual integrability condition
$$\bp(\phi)=-\1/2 [\phi,\phi]\tag 3.7$$
imply that, setting
$$\epsilon (\phi)=(\phi, \phi\times\phi, ...),$$
$$\epsilon(\phi, v)=(v, \phi\times v, \phi\times\phi\times v, ...),$$
the pair
$$\eta (\phi, v):=( \epsilon (\phi), \epsilon (\phi, v))$$
constitute a cocycle for the complex $J_m(\gsh)\otimes\m_s.$
This cocycle is obviously 'morphic' or comultiplicative, hence
gives rise to a ring homomorphism
$$\alpha^\sharp=\alpha^\sharp(\phi,v) : R_m(\gsh)\to S\tag 3.8 $$
lifting the usual Kodaira-Spencer homomorphism $\alpha(\phi):
R_m(\g)\to S$.\par
Conversely, given a homomorphism $\alpha^\sharp$ as above, with
$S$ an arbitrary artin local algebra, clearly we may represent
$\alpha^\sharp$ in the form $\eta (\phi, v)$ as above where $\phi$ and
$v$ satisfy the conditions (3.6) and (3.7). Then in the enveloping
algebra $\U(\h^+)$ we get the identity
$$i(\phi)=-\bp(v)-\1/2 [i(\phi), v].\tag 3.9$$
Plugging the identity (3.9) back into itself we get, recursively,
$$i(\phi)=-\sum\limits_{k=0}^\infty \frac{(-1)^k}{2^k}
\ad (v)^k(\bp(v))\tag 3.10$$
(the sum is finite because $\m_S$ is nilpotent),
from which the identities 3.5 follow formally.
Hence the diagrams 3.3 and 3.4 commute, so we get a $ \gsh-$deformation
lifting $\phi.$\par
In particular, applying this construction to the identity
element of $S=R_m(\gsh)$, thought of as an element of
$\H^0(J_m(\gsh))\otimes\m_S$, we obtain a 'tautological'
$\gsh-$deformation which we denote by
$$u_m^\sharp=(\phi_m,v_m)$$
and we get the following analogue (and extension) of Theorem 3.1:
\proclaim{Theorem 3.2} Let $\gsh$ be an admissible \dg Lie atom
such that $\g$ has central sections. Then for any $\gsh-$
deformation $(\phi,v)$ over an artin local algebra $S$
of exponent $m$ there exist isomorphisms
$$(\phi,v)\simeq \alpha^\sharp(\phi,v)^*(u_m^\sharp)$$
any two of which differ by an element of {\rm{Aut}}$(\phi,v)$.\endproclaim

\vfill\eject
\heading 4. The Hitchin symbol\endheading
Given a nonsingular curve $C$ and a stable vector bundle $E$ on $C$,
Hitchin [Hit] constructed a fundamental map or 'symbol'
$$\Hi_E: H^1(T_C)\to\sym^2(H^1(\sll(E)))\subset \hom(H^0(\sll(E)\otimes K_C),
H^1(\sll(E))).$$
Via the natural identification of $H^1(\sll(E))$ with the tangent space
at $[E]$ to the moduli space $\M$ of bundles with fixed determinant on $C$,
this gives a lifting of the canonical variation map
$$H^1(T_X)\to H^1( \D^1_\M(\theta^k)),$$
where $\theta$ is the theta line bundle on $\M$ and $k$ is
arbitrary,
which is the crucial ingredient in the construction of the flat
connection on the space of sections of $\theta^k$ over $\M$.

Our purpose here is to give a definition of the Hitchin symbol
which is of a 'local' character and which, in particular, avoids
the use of
Serre duality, on which Hitchin's original definition was based.
\subheading {4.1 The definition}
Let
$$\pi :\CC\to S$$
be a family of smooth curves of genus at least 2,
and let $A, B$ be locally free sheaves on $\CC$.
Consider the short exact sequence on $\CC\times_S\CC$:
$$0\to A\bt B\to (A\bt B)(\Delta)\to A\otimes B\otimes T_{\CC /S}\to 0\tag 4.1$$
where $\Delta$ is the diagonal and we have identified
$\O_{\Delta}(\Delta)$ with $ T_{\CC /S}$. This yields a map
of complexes
$$\partial_{A,B}:A\otimes B\otimes T_{\CC /S}\to A\bt B[1]$$
where $A$ and $B$ are identified respectively
with suitable complexes resolving them.
Let
$$\partial_{A,B}^1:R^1\pi_*(A\otimes B\otimes T_{\CC /S})\to R^2\pi_*(A\bt B)$$
be the induced map. Now suppose
given an $\O_{\CC}$-linear 'trace' pairing
$$t:\check{A}\otimes \check{B}\to\O_{\CC} \tag 4.2$$
where $\check{A},\check{B}$ denote the dual locally free sheaves.
This induces
$$\check{t}\otimes\id :T_{\CC /S}\to A\otimes B\otimes T_{\CC /S}.$$
We define the {\it {Hitchin symbol}} associated to
this data to be the composite
$$\hi'_{A,B}=\partial_{A,B}\circ \check{t}\otimes\id: T_{\CC /S}\to
A\bt B[1]$$
or either of the induced maps
$$\Hi'_{A,B/S}=
\partial_{A,B}^1\circ R^1\pi_*(\check{t}\otimes\id ):
R^1\pi_*(T_{\CC /S})\to
R^2\pi_*(A\bt B),$$
$$\Hi'_{A,B}=H^0(\partial_{A,B})\circ H^1(\check{t}\otimes\id ):
H^1(T_{\CC /S})\to H^2(A\bt B).$$
Assuming moreover  that $A=B$ and that $t$ is symmetric, we
get a map
$$\hi_A:T_{\CC /S}\to\lambda^2(A)[1]$$
(note that the shift of 1 exchanges symmetric and alternating
products).
Assuming $\pi_*(A)=0$ we may
 identify
$$R^2\pi_*(A\bt A)= R^1\pi_*(A)^{\otimes 2},$$
and note that $\pi_*(T_{\CC /S})=0$
so we get induced maps
$$\Hi_{A/S}:R^1\pi_*(T_{\CC /S})\to\sym^2(R^1\pi_*(A))$$
i.e. the symmetric component of $Hi'_{A,A/S}$ and likewise for
$$\Hi_A:H^1(T_{\CC /S})\to H^0(\sym^2(R^1\pi_*(A))).$$
\remark{Example 4.1} Let $\g$ be a sheaf of semisimple,
$\O_{\CC}-$locally free Lie
algebras with $\pi_*(\g )=0.$
Then $\g$ is endowed with a nondegenerate trace
pairing
$$t:\sym^2(\g )\to\O_{\CC}$$
which may be used to identify $\g$ and $\g^*$,
whence Hitchin symbols
$$\hi_{\g/S}:T_{\CC /S}\to\lambda^2(\g)[1].$$
$$\Hi_{\g/S}:
R^1\pi_*(T_{\CC /S})\to\H^0(\sym^2(R^1\pi_*(\g))).$$
In particular, if $E$ is a
locally free  $\O_{\CC}-$module we
will abuse notation and denote by $\Hi_{E/S}$
the Hitchin symbol associated to
$\g=\frak {sl} (E).$\par
Specializing further, suppose
$\pi:\CC\to S$ is a given family of
smooth curves of genus $g\geq 3$ endowed with a polarization
of degree $d$, and let $\M\overset {\pi'}\to\to S$
be a locally fine moduli space (cf. \S 6) of stable rank-$r$
bundles of degree $d$ and fixed determinant (= the
polarization) over $\CC/S$; if $(r,d)=1$ we may just take
$$\M=\Cal{SU}^r(\CC/S)\overset {\pi'}\to\to S,$$
the global (fine) moduli space.
Let $E$ be the universal bundle over
$\CC_\M:=\CC\times_S\M$. The we get a Hitchin symbol
(as a bundle map over $\M$):
$$\Hi_{\CC_\M/\M}:\pi^{'*}R^1\pi_* (T_{\CC /S})
\to\sym^2(R^1(\pi\times\pi ')_*(\g))=\sym^2(T_{\M /S})$$
hence, pushing down to $S$ we get a map
$$R^1\pi_* (T_{\CC /S}) \to\pi '_*(\sym^2(T_{\M /S}))
=\pi '_*(\sym^2(T_{\M /S})).$$
This is the map originally defined by Hitchin using Serre duality.

\endremark
The fact that our map and Hitchin's coincide is a consequence of
the following
\proclaim{Proposition 4.2} In the situation
above, the map $\partial^1_{A,B}$ is dual to the restriction map
$$(\pi\times\pi)_*((\check{A}\otimes \Omega_{\CC/S})\bt
(\check{B}\otimes
\Omega_{\CC/S})\to \pi_{\Delta *}(\check{A}\otimes
\check{B}\otimes
\Omega_{\CC/S}^{\otimes 2}).$$\endproclaim
Indeed the Proposition
and the K\"unneth formula imply that
$\Hi_{A,B}$ is dual to the
map $$\pi_*(\check{A}\otimes \Omega_{\CC/S})\otimes
\pi_*(\check{B}\otimes
\Omega_{\CC/S})\to\pi_* ( \Omega_{\CC/S}^{\otimes 2})$$
induced by
multiplication and $t$, which is Hitchin's definition.
\par As for
the Proposition, it follows easily
from relative Serre duality on
the (relative) surface $\CC\times\CC/S$, together with the
following remark
 \proclaim{ Lemma 4.3} Let $D/S\subset X/S$ be an
embedding of a  relative divisor,
with $X,D$ both smooth projective
 over $S$ affine, and let $F$ be a locally free sheaf on $X$.
Then the map
$$H^i(F\otimes\O_D(D))\to H^{i+1}(F)$$
induced by the exact sequence
$$0\to F\to F\otimes\O (D) \to F\otimes\O_D(D)\to 0$$
is dual to the restriction map
$$H^{n-i-1}(F^*\otimes K_X)\to
H^{n-i-1}(F^*\otimes K_X\otimes\O_D).$$
\endproclaim
The Lemma follows easily from any
standard treatment of Serre duality
(e.g. in [Ha]), noting the compatibility of the 'fundamental local
isomorphisms' for $D$ and $X$.\qed
\subheading {4.2 Cohomological formulae}
Our purpose here is to derive
algebraically some cohomological formulae
for Hitchin symbols which
extend and substitute for Hitchin's differential-geometric
calculations in [Hi].
We return to the general situation of (4.1)
above, so $A,B$ are locally free sheaves
on $\CC/S$.
Let $\g$ be a locally free Lie subalgebra of $\gl (A)$. We will
say that $A$ is a $\g-$structure
if we are given a Lie subalgebra
$\gh\subset\D(A)$ which extends
$\g\subset\gl(A)$ (cf. Example 9.1).

Equivalently, as is well known and due to Atiyah,
the Principal Part (or jet) sheaf
$P^1(A)$ can be given as an extension
$$0\to\Omega_{\CC/S}\otimes A\to P^1(A)\to A\to 0$$
with transition cocycle in
$C^1(\g\otimes\Omega_{\CC/S})\subset
C^1(\gl (A) \otimes\Omega_{CC/S})$.
This cocycle then determines the {\it{Atiyah Chern class}}
$$AC(A)\in H^1(\g\otimes\Omega_{\CC/S})$$
(from which the usual Chern classes can be computed); see
also \S 9.
Note that if $\det (A)$ is trivial
then $A$ admits a $\frak{sl}(A)$-structure, where
 $\frak{sl}(A)=\{$ endomorphisms of $A$ acting trivially on
$\det (A)\}$; also, for any $\O_{\CC}$-Lie algebra
bundle $\g$, $\g$ itself
admits a $\g$-structure, via the adjoint representation.
\proclaim{Lemma 4.4} For arbitrary locally free sheaves $A,B$
on $\CC$, there
is a natural $\O_{\CC}-$ isomorphism
$$p_{1*}(A\bt B\otimes\O_{2\Delta}(\Delta))\simeq
A\otimes P^1(B\otimes T_{\CC/S})$$
where $p_i$ denote the coordinate projections.\endproclaim
\demo{proof} To begin with, we have
$$p_{1*}(A\bt B\otimes\O_{2\Delta}(\Delta))=
p_{1*}(p_1^*A\otimes (p_2^*B\otimes\O_{2\Delta}(\Delta)))$$
$$\simeq A\otimes p_{1*}(p_2^*B\otimes
\O_{2\Delta}(\Delta)).$$
Hence we may assume $A=\O_{\CC}.$
Next, note the natural map
$$p_{2*}(p_2^*B\otimes\O(\Delta)|_{\Delta})
\to p_{2*}(p_2^*B\otimes
\O(\Delta)\otimes\O_{\Delta})=B\otimes
T_{\CC/S}$$
where $|_{\Delta}$ denotes topological restriction.
This gives rise to a map
$$p_2^*p_{2*}(p_2^*B\otimes\O(\Delta)|_{\Delta})\to
p_2^*(B\otimes T_{\CC/S})$$
But since $p_2:\Delta\to\CC$ is a homeomorphism we may identify
$$p_2^*p_{2*}(p_2^*B\otimes\O(\Delta)|_{\Delta})=
p_2^*B\otimes\O(\Delta)|_{\Delta}$$
so we get a map
$$p_2^*B\otimes\O(\Delta)|_{\Delta}\to  p_2^*(B\otimes T_{\CC/S})$$
which induces a map
$$\rho: p_2^*B\otimes\O(\Delta)\otimes\O_{2\Delta}\to
p_2^*(B\otimes T_{\CC/S})\otimes\O_{2\Delta}=
P^1(B\otimes T_{\CC/S}).$$
Now both sides admit (2-step) filtrations with the same graded
pieces, and it is easy to see locally that $\rho$ is compatible
with these filtrations and induces the identity on the
gradeds, hence  is an isomorphism.\qed
\enddemo
Now it is well known and easy to prove that, via the natural
inclusions
$$\endf (F), \endf(B) \subseteq \endf (F\otimes B),$$
$$a\mapsto a\otimes\id, b\mapsto \id\otimes b,$$
the holonomy algebra of $F\otimes B$ is generated by those of
$F$ and $B$ and
we have
$$AC(F\otimes B)=\rk (F)AC(B)+\rk (B)AC(F).$$
Consequently, we have
\proclaim{Corollary 4.5} The map
$$A\otimes B\otimes T_{\CC /S}\to A\otimes B[1]$$
induced by $\partial_{A,B}$ and multiplication is given
by cup product
with
$$-\rk (B)c_1(K_{\CC /S})+AC(B)$$
where $c_1(K_{\CC /S})\in H^1\pi_*(\Omega_{\CC /S})$
is the relative canonical class.\endproclaim

\proclaim{Corollary 4.6} In the presence of s symmetric
trace map as in display (4.2), the map
$$T_{\CC /S}\to A\otimes B[1]$$
induced by $\Hi_{A,B}'$ and multiplication is given by
$$-\rk (B)\check{t}\circ c_1(K_{\CC /S})
+AC(B)\circ (\check{t}\otimes\id_{ T_{\CC /S}}).$$
\endproclaim
\proclaim{Corollary 4.7} If $\beta:A\otimes A\to C$
is any skew-symmetric
pairing, then the map
$$ T_{\CC /S}\to C[1]$$
induced by $\Hi_{A}$ and $\beta$ is given by
$$\beta\circ AC(B)\circ (\check{t}\otimes\id_{ T_{\CC /S}}).$$
\endproclaim
\demo{proof} It suffices to note that $\check{t}$
lands in the symmetric part
of $A\otimes A$, hence
$\check{t} \circ c_1(K_{\CC /S})$ is mapped to zero by $\beta$,
so the previous corollary yields
the result.\qed\enddemo
\proclaim{Proposition 4.8} Let $E$ be a locally free sheaf
on $\CC/S$ and let
$$\ell\in \g^1\otimes\Omega_{\CC /S},\ \g=\frak{sl}(E).$$
be a representative of the traceless part of the Atiyah-Chern class
$AC(E).$
 Then the map
$$ T_{\CC /S}\to \g[1]$$
induced by $\Hi_E$ and the bracket on $\g$ is
given by $2\rk(E)\ell.$\endproclaim
\demo{proof} Let $(e_i)$ be a local frame for $E$ and $(e_j^*)$
the dual frame. Then modulo scalars
$\ell$ can be written in the form
$$\ell_{\alpha\beta} =
\sum\limits_{i,j}\phi_{\alpha\beta i j}e_i^*\otimes e_j$$
for certain local sections $\phi_{\alpha\beta i j}$ of
$\Omega_{\CC /S}$
and $AC(E\otimes E^*)$ can be written as
$$\widetilde{\ell_{\alpha\beta}}=
-\sum\limits_{j,k,l}\phi_{\alpha\beta k j}e_j\otimes
e_l^*\otimes e_k^*\otimes e_l+
\sum\limits_{i,k,l}\phi_{\alpha\beta i l}e_k\otimes e_i^*
\otimes e_k^*\otimes e_l.$$
Now the dual $tr^*$ of the trace pairing can be expressed as
$$tr^*=\id\otimes\sum\limits_{l,k}e_l\otimes e_k^*
\otimes e_k\otimes e_l^*.$$
Therefore by Corollary 4.2.3 the Hitchin symbol
$Hi'_{E\otimes E^*, E\otimes E^*}$,
except for a symmetric part which is killed by the bracket,
is given by
$$\sum\limits_{k l j} \phi_{\alpha\beta k j}e_l\otimes e_k^*
\otimes e_j\otimes e_l^*-
\sum\limits_{k l i} \phi_{\alpha\beta i l}e_l\otimes e_k^*
\otimes e_k\otimes e_i^*.$$
Now the result of applying bracket to this is that of
contracting the
middle two factors minus that of contracting the
outer two. Now contracting
the inner (resp. outer) factors on the first
(resp. second) sum yields
a multiple of the identity which may be ignored. The rest yields
$$\rk (E)\sum\limits_{k  j} \phi_{\alpha\beta k j} e_k^*\otimes e_j
+\rk (E)\sum\limits_{ l i} \phi_{\alpha\beta i l}e_l\otimes e_i^*$$
$$=2\rk (E)\ell_{\alpha\beta}.\qed$$
\
\enddemo
\vfill\eject
\heading 5. Moduli modules revisited\endheading

Let $(\g,E)$ be an admissible pair with
$H^0(\g)=0$ on a Hausdorff
space $X$. In [ Rcid ] we constructed
the universal deformation ring
$$\hat{R}(\g)=\lim\limits_{\leftarrow} R_m(\g),$$
$$R_m(\g)=\C\oplus \H^0(J_m(\g))^*$$
where $J_m(\g)$ is the Jacobi complex of $\g$, as well as the
flat $R_m(\g)$-module $E_m=M_m(\g,E)$ that is the universal
$\g$-deformation of $E$, and whose
cohomology groups may be called the
'moduli modules' associated to $E$. Our purpose here is to
revisit that construction from a slightly different
viewpoint that seems more convenient for applications, such
as the construction of Lie brackets on moduli.\par
As in [Rcid] we let $(\g^.,\delta )$ be a soft dgla
resolution of $\g$ and
$(E^.,\partial )$ be a soft resolution of $E$ that is a graded
$\g^.$-module. Then the standard (Jacobi) complex
$J_m(\g^.)=:\J_m^.$ has
terms which may be decomposed as
$$\lambda^i(\g^.)=\bigoplus\limits_j \g^{-i}_j $$
where each $\g^{-i}_j $ has total degree $j-i$
and is a sum of terms of
the form
$$(\lambda^{a_1}\g^{b_1}\bt\cdots)\bt
(\sigma^{c_1}\g^{d_1}\bt\cdots)$$
with $b_k$ even, $d_k$ odd and
$$\sum a_k+\sum c_k=i, \sum a_kb_k +\sum c_kd_k=j.$$
 Thus $\H^.(J_m(\g))$ is $H^.$ of a complex
$\Gamma J^._m$ with
$$\Gamma J^r_m =\bigoplus\limits_{i=1}^m
\Gamma (\g^{-i}_{r+i}) .$$
It is convenient to augment $\Gamma J^r_m$ by adding the term
$\g^0_0=\C$ in degree 0.\par
Note that $\Gamma J^r_m$ depends only on
the \dgla  $\ g^.=\Gamma (\g^.)$; moreover
the quasi-isomorphism class of $\Gamma J^r_m$
depends only on the dgla quasi-isomorphism class of  $\g^.$.
 Also,
it is worth noting at this point that
the differential $\delta$ on $\J_m^.$ and
$\Gamma J^._m$ is a 'graded coderivation', in the sense of
commutativity of the following diagram
$$\matrix \J_m^.&\overset\delta\to\to&\J_m^.\\
\downarrow&&\downarrow\\
\sigma^2\J_{m-1}^.&
\overset{\sigma^2(\delta)}\to\longrightarrow
&\sigma^2\J_{m-1}
\endmatrix $$
where the vertical arrows are the comultiplication maps
and $\sigma^2(\delta)$ is induced by $\delta$ by functoriality,
i.e. is the map given by extending $\delta$ 'as a derivation',
given by the rule $$ab\mapsto \delta(a)b\pm a\delta(b).$$
Thus $\J_m^.$ and
$\Gamma J^._m$ possess are 'differential graded
OS structures . It follows that
the same is true not only
of their $H^0$ but of the entire cohomology
$$H^.(\Gamma J^._m)=\bigoplus H^i(\Gamma J^._m)=
\bigoplus\H^i(J_m(\g )).$$
Now to get an algebra out of this one may either dualize the
cohomology, as was done in [Rcid] or,
what essentially amounts to the same
thing but seems more convenient, one may dualize the complex
and then take cohomology. Thus let
$\Gamma^*J^._m$ be the complex dual to
$\Gamma J^._m$, with
$$\Gamma^*J^r_m=(\Gamma J^{-r}_m)^*$$
(vector space dual) and dual differentials. Thus
$$ \Gamma^*J^r_m=\bigoplus
\limits_{i=1}^m\Gamma (\g^{-i}_{-r+i})^*,$$
$$\Gamma(\g^{-i}_j)^*=
\bigoplus [\bigwedge^{a_1}\Gamma(\g^{b_1})\otimes
\cdots]\otimes[\sym^{c_1}\Gamma(\g^{d_1})\otimes\cdots],$$
with sum over all nonnegative
 with $b_k$ even, $d_k$ odd and
$$\sum a_k+\sum c_k=i, \sum a_kb_k +\sum c_kd_k=j.$$
 Clearly
$$H^i(\Gamma^*J^._m)=H^{-i}(\Gamma J^._m)^*.$$
Then $\C\oplus\Gamma^*J^._m$ is in fact a
differential graded- commutative
associative algebra (graded commutative means two homogeneous
elements commute unless they are both odd, in which case they
anticommute). Indeed, since the OS or comultiplicative
structure on $J$ was derived from exterior comultiplication
on $\lambda^.(\g^.)$,
clearly the multiplication on
$\C\oplus\Gamma^*J^._m$ is induced by
graded exterior multiplication, hence
is obtained by tensoring together the
various exterior products on the
$\bigwedge^{a_k}\Gamma(\g^{b_k}), b_k$
even and symmetric products on
the $\sym^{c_k}\Gamma(\g^{d_k}), d_k$ odd.
Thus the total cohomology
$$\tilde{R}_m(\g):= H^.(\Gamma^*J^._m)=
(H^.(\Gamma J^._m)^*)$$
is a local graded artin algebra, and in
particular the degree-0 piece
$$R_m(\g)= H^0(\Gamma^*J^._m)= (H^0(\Gamma J^._m)^*)$$
is the $m$th universal deformation ring mentioned above.\par
Next, we revisit the construction of the $m$-universal
deformation $M_m(\g,E)$ of a given $\g-$module $E$.
The proof of ([Ruvhs], Thm 3.1)shows that
the MOS structure
$V^m(E)$- which in turn determines $M_m(\g,E)$ purely formally-
is a cohomology sheaf $\SH^0$ of a (multiple) complex
$(K^._m,d^.)=(K^._m(\g,E),d^.)$ whose part in
total degree $r$ is
$$K^r_m=\bigoplus\limits_{{i,j, j\geq-m}}
\Gamma(\g^j_{r-j-i})\otimes E^i
= \bigoplus\limits_{{i}}\Gamma J^{r-i}_m\otimes E^i$$
(note there is a misprint in the corresponding formula
in ([Ruvhs], p.430, l.-5).
'Transposing' this, we define a complex
$$L^.=^tK^._m(\g,E)$$ of sheaves with
$$L^r=\bigoplus\limits_{i,j}\Gamma
(\g^{-j}_{i+j-r})^*\otimes E^i
 =\bigoplus\limits_{i}\Gamma^*J^{r-i}_m\otimes E^i$$
with differentials the 'transpose' of those of $K^.$, where the
transpose of a map
$$d:A\otimes E\to B\otimes E'$$
with $A,B$ finite-dimensional vector spaces,
is the map
$$^td:B^*\otimes E\to A^*\otimes E'$$
defined by the rule
$$<^td(b^*\otimes e),a>=<b^*,d(a\otimes e)>$$
in which $<,>$ refers to the natural pairings
$$A^*\otimes E'\times A\to E', B^*\times B\otimes E'\to E'.$$
Since $\Gamma J_m$ has finite-dimensional cohomology we
may 'approximate' it by a finite-dimensional subcomplex quasi-
isomorphic to it (this remark will be
used frequently in the sequel).
Since the comultiplicative structure on $K^.$  as defined
in [Ruvhs] coincides with the evident structure induced by
comultiplication
on the $\g$ factor, clearly $L^.$
has a natural structure of a sheaf of differential graded
$\C\oplus\Gamma^*J^._m$-modules,
and consequently the total cohomology
$\SH^.(L^.)$ is a sheaf of graded
$\tilde{R}_m(\g)$-modules and in
particular $\SH^0(L^.)$ is a sheaf of $R_m(\g)$-modules.
Note also that for any $\g-$modules $E_1, E_2$,
the multiplicative structure on $\Gamma^*J^._m$
gives rise to a natural pairing
$$^tK^._m(\g,E_1)\times ^tK^._m(\g,E_2)\to
^tK^._m(\g,E_1\otimes E_2).\tag 5.1$$
Note also that $\Gamma L^.$ depends only on
$\Gamma E^.$ (as \dg\ $g^.)-$
module), and we may similarly associate a
complex, still denoted
$^tK^._m(g^.,F^.)$, to any  \dg\ $g^.-$module $F^.$.
\proclaim{Theorem 5.1} We have natural isomorphisms
$$E_m\simeq M_m(\g,E)\simeq \SH^0(^tK^._m(\g,E)).$$\endproclaim
\demo{proof} The first isomorphism is just
[Ruvhs], Thm 3.1, but it is
worth observing that the proof given there involves an implicit
spectral sequence argument, and may be shortened by making this
argument explicit. Thus, write $K^.$ as a double complex
$$ K^{i,j}=\bigoplus\limits_{k}
\Gamma(\g^k_{j-k})\otimes E^i=\Gamma J_m^j\otimes E^i. $$
Now because $H^0(\g)=0$, the map
$$\delta:\Gamma(\g^0)\to\Gamma (\g^1)$$
is injective, so choosing a complement to its image we obtain
a quasi-isomorphic complex in strictly positive degrees, and
it will be convenient to replace all resulting complexes by ones
formed with this modified complex. This in particular ensures
that $K^{i,j}=0$ for $j<0.$
Note that the vertical differentials
$$K^{i,j}\to K^{i,j+1}$$
are of the form $\delta\otimes id$ where
$\delta$ is a differential of
$\Gamma J^._m$.
Consequently, the first spectral sequence of the double
complex has an $E_1$ term
$$\matrix E_1^{p,q}&=&(\C\oplus H^0(\Gamma J^._m))
\otimes E^p, q=0\\
                   &=& H^q(\Gamma J^._m)\otimes E^p, q>0.
\endmatrix$$
Our assumption that $H^0(\g)=0$ easily implies that
$ H^q(\Gamma J^._m)=0$ for $q<0$, so
this is a first-quadrant spectral sequence and consequently
$$H^0(K^.)=\ker (E_1^{0,0}\to E_1^{0,1})=
\ker (V_0^m\otimes E^0\to V_0^m\otimes E^1),$$
and identifying the map and applying a suitable
functor gives the result.\par
For the second isomorphism we argue analogously,
considering the double complex
$$L^{i,j}=\bigoplus\limits_{k}
\Gamma(\g^k_{-j-k})^*\otimes E^i=
\Gamma^* J_m^j\otimes E^i. $$
As above, this vanishes for $j>0.$
We get a spectral sequence whose $E_1$ term
may be identified as
$$E_1^{p,q}=\tilde{R}_m(\g)^q\otimes E^p$$
where $\tilde{R}_m(\g)^q$ denotes the part in degree $q$
i.e. $H^q(\Gamma^*J^._m)$ for $q>0$ and
$\C\oplus H^0(\Gamma^*J^._m)$ for $q=0$ .
Thanks again to our hypothesis that $H^0(\g)=0$,
this vanishes for all
$q>0$, so in this case we have a fourth-quadrant
spectral sequence.
Now note that if we view
$(\bigoplus\limits_{q}E_1^{p,q})$ as a complex
indexed by $p$ only, the differentials are
$\tilde{R}_m(\g)$-linear,
so as such this is a finite complex of free
$\tilde{R}_m(\g)$-modules. Now
we use the following fairly standard fact.
\proclaim{Lemma 5.2} Let $A^.$ be a  complex
of flat modules
(not necessarily of finite type) over a local artin algebra $S$
with residue field $k$.
Suppose $A^.\otimes_S k$ is exact in positive
degrees. Then so is $A^.$\endproclaim
\demo{proof} Use induction on the length of $S$. Let $I<S$ be
a nonzero 'socle' ideal, with $I.\m_S=0$.
By flatness we have a short exact sequence of
complexes
$$0\to IA^.\to A^.\to A^.\otimes (S/I)\to 0.$$
As $IA^.$ is isomorphic to a direct sum of
copies of $A^.\otimes_S k$,
it is exact in positive degrees; by induction, the same is true of
$A^.\otimes (S/I)\to 0$.
By the long cohomology sequence, it follows that
the middle is also exact in positive degrees,
proving the Lemma.\qed\ \
(The Lemma is perhaps more familiar in the case
where the $A^.$ are of finite type (hence free)
and zero for $\cdot>>0$, and $S$ is not necessarily
artinian.)\enddemo
  For our complex above tensoring
with the residue field just yields the original complex $E^.$
which is of course exact in positive degrees,
so Lemma 5.2 applies. We conclude that
$E_1^{p,q}$ is exact at all terms with $p>0$
and in particular we have
$$H^0(L^.)=\ker (E_1^{0,0}\to E_1^{0,1})=
\ker (R_m(\g)\otimes E^0\to R_m(\g)\otimes E^1),$$
and again the map may be easily identified
as the one yielding
$E_m$.\qed

\enddemo
We will now formalize some constructions which
occurred in the foregoing proof.
For a double complex $K^{..}$ we will denote by
$K^{..}_{\lceil j }$ its $j$th lower truncation, which
is the double complex defined by
$$K^{h,i}_{\lceil j }=\{\ \   \matrix 0,& i>j,\\
\ker (K^{h,j}\to K^{h,j+1}),& i=j\\
K^{h,i}, &i<j.\endmatrix$$
Similarly,  we will denote by
$K^{..}_{\lfloor j }$ the $j$th upper truncation, which
is the double complex defined by
$$K^{h,i}_{\lfloor j }=\{\ \   \matrix 0,& i<j,\\
K^{h,j}/K^{h,j-1},& i=j\\
K^{h,i}, &i>j.\endmatrix$$
Motivated by the foregoing proof, we set
$$\Lt_m(\g,E)=\ ^t(K_m^.(\g,E)) _{\lceil 0}\tag 5.2$$
which, as we have seen, is a double complex in nonpositive
vertical degrees (indeed in the fourth quadrant)
  quasi-isomorphic to $^t(K_m^.(\g,E))$ itself.
Also set
$$L_m(\g,E)=\Lt_m(\g,E)_{\lfloor 0}\tag 5.3$$
which is thus to be considered as a simple (horizontal)
complex.
\proclaim{Corollary 5.3} Assumptions as in 5.1, we have
\item{(i)} $$\SH^i(L_m^.(\g,E))=\{ \matrix 0 &i\neq 0\\
M_m(\g,E)& i=0\endmatrix$$
\item{(ii)}$$H^i(X,M_m(\g,E))=H^i(\Gamma(L_m^.(\g,E)))$$
\endproclaim\qed
\par
Note that the above constructions make sense for any  \dgla
$g^.$ and \dg $g^.$-module $G^.$ and yield (double) complexes
$$K_m(g^.,G^.):\ \ \ \  K^{i,j}=J_m^j(g^.)\otimes G^i$$
$$^tK_m(g^.,G^.):\ \ \ \  ^tK^{i,j}=^*J_m^j(g^.)\otimes G^i.$$
Likewise for $\Lt_m, L_m$.
Moreover, clearly the quasi-isomorphism classes of these
complexes depends only on
that of $G^.$ as $g^.$-module. We collect some of the
simple properties of these constructions
in the following
\proclaim{Lemma 5.4} In the above situation, assume
$H^{\leq 0}(g^.)=0$.Then
\item {(i)} If $G^.$ is acyclic in negative degrees, then so is
$K_m(g^.,G^.)$(i.e. it is acyclic in negative total degrees).
\item{(ii)} If $G^.$ is acyclic in positive degrees,
then so is $^tK_m(g^.,G^.)$.
\item{(iii)}If $G^.$ is acyclic in negative degrees,
then $\Lt_m(g^.,G^.)$ is a 4th- quadrant bicomplex and
$L_m(g^.,G^.)$ is acyclic in negative degrres..
\item {(iv)} The complex $^*G^.$ dual to $G^.$ is also
a $g^.$-module, and
we have
$$ ^* (K_m(g^.,G^.))=^tK_m(g^.,^*G^.).$$\endproclaim\qed
\proclaim{Corollary 5.5} In the situation of Corollary 5.3,
assume moreover that
for some $i,$
$$H^j(E)=0, \forall j\neq i.$$
Then we have (where * denotes vector space dual)\ss
\item{(i)}\centerline{
$ H^i(K_m(\g,E))^*= H^{-i}(^tK_m(\Gamma(\g^.),\Gamma(E^.)^*)
=H^{-i}(L_m(\Gamma(\g^.),\Gamma(E^.)^*)),$}\ss
\item{(ii)}\centerline{
$H^i(M_m(\g,E))^*= H^{-i}(K_m(\Gamma(\g^.),\Gamma(E^.)^*)).$}
\endproclaim
\demo{proof} (i) The first equality is just the
fact that cohomology commutes with dualizing. The second
follows by applying Lemma 5.4 to $\Gamma (E^.)[-i]$ which is
acyclic except in degree 0 and consequently
$H^{-i}(^tK_m(\Gamma(\g^.),\Gamma(E^.)$ only involves
$H^0(\Gamma J_m)$. (ii) follows from Corollary 5.3.\qed\enddemo
A sheaf (or complex) satisfying the condition of Corollary 5.5
will be said to be {\it{equicyclic of degree $i$}}
or $i$-{\it equicyclic}.
Next, note that a $g^.-$bilinear pairing of \dg\  $g^.-$modules
$$G_1\times G_2\to G_3$$
gives rise to a pairing of double complexes
$$\Lt_m(G_1)\times \Lt_m(G_2)\to \Lt_m(G_3)\tag 5.4$$
hence also
$$L_m(G_1)\times L_m(G_2)\to L_m(G_3)\tag 5.5$$
Clearly these pairings are compatible with the $^*J_m(g^.)$-module
structure on these complexes, so we get an $R_m(g^.)$-
linear pairing
$$H^i(L_m(G_1))\times H^j(L_m(G_2))\to H^{i+j}(L_m(G_3)).$$
In particular, for any \dg $g^.$-module $G^.$, using the
natural $g^.$-linear
pairing
$$G^.\times^*G^.\to \C$$
(where $\C$ is endowed with the trivial $g^.$-action), we obtain
an $R_m(g^.)-$linear pairing
$$H^i (L_m(G^.))\times H^{-i}(L_m(^*G^.))\to
H^{0}(L_m(\C))=R_m(g^.).\tag 5.5$$
Now suppose moreover that we have $G^.$is  $i-$equicyclic.
Then clearly $H^i (L_m(G^.))$ is a free $R_m(g^.)-$module
(as obstructions lie in $H^{i+1}(G^.)=0$ ) and similarly
for $H^{-i}(L_m(^*G^.) )$. Since the pairing (5.5) yields the
natural perfect
pairing of $H^i(G^.)$ and $H^{-i}(^*G^.)$ modulo the
maximal ideal of
$R_m(g^.)$, it is likewise perfect. Thus
\proclaim{Corollary 5.6} In the above situation, if $G^.$ is\
\ $i-$equicyclic, then
$$H^i (L_m(G^.))\ {\text{ and}}\ \ H^{-i}(L_m(^*G^.))$$
are free $R_m(g^.)-$modules
naturally dual to each other.\endproclaim
\proclaim{Corollary 5.7}
  In the situation of Corollary 5.5, $$H^i(K_m(\g,E))$$
  is the $R_m(\g)$-module dual to
the free module $$H^{-i}(K_m(\Gamma(\g^.),\Gamma^*(E^.))).
$$\endproclaim
\demo{proof} This follows from Corollary 5.6 and the fact that
$H^i(K_m(\g,E))$ and Hom$(H^{-i}(L_m(\Gamma^*(E))), \C)$ coincide as
$R_m$-modules.\qed\enddemo
\proclaim{Lemma 5.8} For any $g^.-$module $G^.$,
\item{(i)} $K_m(G)$ and $L_m(G)$ are $g-$modules;
\item{(ii)}there is a natural inclusion $L_{m+k}(G)\to L_k(L_m(G)).$
\endproclaim
\demo{proof}(i) The point is that the natural $g-$
action on the components
of $K(G)$ and $L(G)$ (suppressing $m$ for convenience)
commutes with
the differentials. Firstly for $K$ and for its first
differential
$$K^{-1}(G)=g\otimes G\to K^0(G)=G,$$ this commutativity is verified
by the fact that
$$<v,<w,a>>=<[v,w],a>+<w,<v,a>>, \forall v,w\in g, a\in G,$$
which means that the following diagram commutes
$$\matrix g\times g\otimes G&\overset{\id\times\delta}\to
\to& g\times G\\
\downarrow&&\downarrow\\
g\otimes G&\overset{\delta}\to\to&g\endmatrix$$
(vertical arrows given by the action; NB $g$ acts both on
$g$ (adjoint
action) and $G$).\par
Next, the case of an arbitrary differential of $K$ follows
by noting the inclusion $K^{-i}(G)\subset
K^{-1}(\bigwedge\limits
^{i-1}g\otimes G), \forall i$, which makes
the following diagram commute
$$\matrix K^{-i}(G)&\overset\delta\to\to&K^{-i+1}(G)\\
\downarrow&&\parallel\\
K^{-1}(\bigwedge\limits^{i-1}g\otimes G)&
\overset\delta\to\to&
K^0(\bigwedge\limits^{i-1}g\otimes G).\endmatrix$$\par
For the case of the $L$ complex, again it suffices to prove
commutativity of the action with the first differential,
i.e. commutativity of the diagram
$$\matrix g\times G&\overset{id\times\cob}\to\to&
g\times ^*g\otimes G\\
\downarrow&&\downarrow\\
G&\overset{\cob}\to\to&^*g\otimes G,\endmatrix$$
which amounts to
$$\cob(<v,a>)=<v,\cob(a)>,\forall a\in G,v\in g.\tag 5.6$$
This is verified by the following computation.
Pick $y\in g$ and write
$$\cob(<v,a>)=\sum w_i^*\otimes b_i.$$
Then (assuming $a,v,y$ are all even)
$$<y.\cob(<v,a>)>=<y,<v,a>>=<[y,v],a>+<v,<y,a>>$$
$$=<[y,v],a>+<v,<y.\cob(a)>$$$$=<[y,v],a>+<[v,y].\cob(a)>+
<y.<v,\cob(a)>>$$
$$=<[y,v],a>+<[v,y],a>+<y.<v,\cob(a)>>=<y.<v,\cob(a)>>.$$
Similar computations can be done for other parities.
Thus (5.6) holds, as claimed.\par
(ii)\ Note that $L_k(L_m(G))$ is naturally a double complex
with vertical differentials those coming from $L_m(G)$. Then
each term of the associated total complex is a sum of
copies of $\bigwedge\limits^{i}(^*g)\otimes G$ and naturally
contains $\bigwedge\limits^{i} (^*g)\otimes G$ itself, embedded
diagonally. It is easy to check that this yields a map
of complexes $L_{m+k}(G)\to L_k(L_m(G)).$\qed
\enddemo
\remark{Remark} The latter inclusion is analogous,
and closely related to, the natural map on jet or principal
parts modules
$$P^{m+k}(M)\to P^k(P^m(M))$$
for any module $M$ (over a commutative ring).\endremark
\vfill\eject

\heading 6. Tangent Algebra\endheading
The purpose of this section is to construct the
tangent (or  derivation)
{\it {Lie algebra}} of vector fields on a moduli space
$\M$, together with its natural representation on the (formal)
functions on $\M$. More specifically, we will say that a locally
$\C-$ringed topological space
$\M$ is a {\it{locally  fine moduli space}} if there exists an
$\O_X$- Lie algebra
$\tilde{\g}$ and a sheaf $\E$ of $\O_X$ and
$\tilde{\g}$-modules on $\M\times X$, with $\tilde{\g}$
acting $\O_X-$linearly,
such that for each point $[E]\in\M,$
we may identify
$$ E=\E|{[E]\times X}:=\E\otimes\C (E)$$
($\C(E)=$ residue field of $\M$ at $[E]$), and the
formal completion
$\hat{\E}$ of $\E$ along $[E]\times X$ is isomorphic
to the universal formal
$\g_E$-deformation of $E$, where
$\g_E=\tilde{\g}\otimes\C (E)$, as constructed
in [Rcid] and above, so that for each $m$, or at least a
cofinal set of $m$'s, we have  (compatible) isomorphisms
$$\E\otimes (\O_{\M}/\m_{[E]}^{m+1})\simeq
M_m(\g_E,E),\tag 6.1$$
$$\O_{\M}/\m_{[E]}^{m+1}\simeq R_m(\g_E).\tag 6.2$$
We do {\it{not}} assume points of $\M$
correspond bijectively with
'equivalence' classes of objects $[E]$
(which we don't even define)--
when a fine moduli space $\frak M$ does exist, our assumptions
imply that the natural classifying  map
$\M\to\frak M$ is \'etale.
Of course by definition the above properties
essentially depend on $\tilde{\g}$ only and not on
the particular $\tilde{\g}$ -module  $\E$. Then
the
tangent sheaf
$$T_{\M}\simeq R^1p_{1*}(\tilde{\g})\tag 6.3$$
(isomorphism as $\O_{\M}$-modules).
Now fix a point $[E]\in \M$ and set $\g=\g_E.$
Viewing $\tilde{\g}$ as a module
over itself via the adjoint representation, we get an isomorphism
of the jet or principal part space
$$P^m(T_{\M})\otimes\C (E)\simeq M_m(\g,\g).\tag 6.4$$
Now the Lie bracket on $T_{\M}$ is a first-order differential
operator in each
argument, hence yields an $\O_{\M}$-linear 'bracket' pairing
$$B_m: P^m(T_{\M})\times P^m(T_{\M})\to P^{m-1}(T_{\M}).\tag 6.5$$
Likewise, the action of $T_{\M}$ on $\O_{\M}$ yields an
'action' pairing
$$A_m: P^m(T_{\M})\times \O_{\M}/\m_{[E]}^{m+1}\to
\O_{\M}/\m_{[E]}^{m}\tag 6.6.$$
The problem is to identify the pairings (6.5), (6.6)
in terms of the
identifications  (6.2) and (6.4). We shall proceed to
define some pairings on
complexes that will yield this. \par
First,
it was already observed above that the dgla $\Gamma(\g^.)$
is quasi-isomorphic
to a sub-dgla of itself in strictly positive degrees.
More canonically, we may set
$$\Gamma^i=\{\matrix 0,&\ \ i\leq 0\\
        \g^1/\delta (\g^0),&\ \ i=1\\
        \Gamma (\g^i)/[\delta\Gamma (\g^0),\Gamma (\g^{i-1})],
        & \ \ i>1.\endmatrix$$
 Then $\Gamma^.$ is a canonical quasi-isomorphic quotient
 dgla of $\Gamma(\g^.)$
in positive degrees. Although a given $\g-$module
$E$ may not give rise to a
$\Gamma^.-$module, still for the purposes of this
section we may as well
replace $\Gamma (\g^.)$ by $\Gamma^.$ and assume it exists
only in positive degrees. Let us also set, for convenience
$$g^.=\Gamma(\g^.).$$
\par
We now begin constructing the action pairing.
Note that the complex $^*g^.=\Gamma(\g^.)^*$ is naturally
a graded module over the dgla $g^.$ known as the
{\it{coadjoint representation}}, via the rule
$$<<a,b^*>,b>=<b^*,[a,b]>, a,b,\in\Gamma(\g), b^*\in\Gamma^*(\g).$$
Hence we get a complex which we will write as $^tK_m(\g,\g^*)$
or $^tK_m(g^.,^*g^.).$
To abbreviate, we will also write $\Gamma(^tK_m(\g,E))_+$
as $\Lt_m(E)$
and $^tK_m(\g,\g^*)_+$ as $\Lt_m(\g^*)$, and we will view them as
double complexes
in nonpositive vertical degrees (in the latter case,
nonpositive horizontal degrees as well).
One can check easily that the duality pairing
$$g^.\otimes ^*g^.\to\C,$$
viewed as a map between $\g^.$-modules (where $\C$ has the
trivial
action), is $g^.$-linear, hence gives rise to a pairing
of double complexes
(preserving total bidegree)
$$\Lt_m(\g)\times \Lt_m(\g^*)\to \Lt_m(\g\otimes\g^*)
\to \Lt_m(\C)=\Gamma^*J_m\tag 6.7$$
(where the RHS is viewed as a double complex in bidegrees
$(0,\cdot\leq 0)$).
Next, note the natural map
$$\Gamma^*J_{m+1}\to \Lt_m(\g^*)[-1],$$
where the shift is taken vertically. This map comes about
by writing symbolically
$$\Gamma^*J_{m+1}: \C\overset{0}\to\to
^*g\to\bigwedge\limits^{2}\ ^*g\to\ldots$$
$$\Lt_m(\g^*)[-1]: ^*g\to^*g\otimes^*g\to\ldots,$$
and mapping
$\C$ to 0 and $\bigwedge\limits^{i}\ ^*g\to
\bigwedge\limits^{i-1}\
^*g\otimes^*g$
in the standard way.
According to our conventions, this map only preserves total
degree; it induces
$$\Lt_{m+1}(\C)=\Gamma^*J_{m+1}\to \Lt_m(\g^*)[-1].$$
Combining the latter with the pairing (6.7), we get a pairing
$$<.>:\Lt_m(\g)\times \Lt_{m+1}(\C)\to \Lt_m(\C)[-1].$$
Now this map takes an element of bi-bidegree $((a_1,a_2),(0,b))$
to a sum of elements of bidegrees $(a_1+b_1=0,a_2+b_2+1)$, where
$b_1+b_2=b$ and $(b_1,b_2)$ is the bidegree of an element in
$L_m(\g^*).$
 Since $ a_2,b_2\leq 0$, it follows that $a_2+b_2+1<1$
if either $a_2<0$ or $b<0$. Therefore there is an induced map
$$<.>:L_m(\g)\times L_{m+1}(\C)\to L_m(\C)[-1],\tag 6.8$$

whence a pairing on cohomology
$$H^1(L_m(\g))\times H^0(L_{m+1}(\C))\to H^0(L_m(\C)).$$
Note that $H^0(L_m(\C))=R_m(\g).$ We set
$$\Theta_m(\g)=H^1(L_m(\g)).$$
By Corollary 5.3, this group coincides with $H^1(M_m(\g,\g))$, i.e.
the $m$-th principal part of the tangent sheaf to moduli.
Thus we have defined a pairing (action pairing)
$$<.>:\Theta_m(\g)\times R_{m+1}(\g)\to R_m(\g).\tag 6.9$$
Now it is easy to see that the pairing
$$L_m(\C)\times L_m(\g)\to L_m(\g)$$
(which comes from the 'product of L's maps
to L of product' rule (5.4)) induces
$$R_m(\g)\times \Theta_m(\g)\to \Theta_m(\g)$$
which endows $\Theta_m(\g)$ with an $R_m(\g)$-module structure.
\par
Next, we undertake to define a pairing on $\Theta_m(\g)$
that will yield
the Lie bracket. For this consider the complex
$\Lt_m(\sym^2g^.)$ which
may be written in the form
$$\sym^2g^.\to (\Gamma^*J^{(1)}_m)_{[-1}\otimes\sym^2\g^.
\to(\Gamma^*J^{(2)}_m)_{[-2}\otimes\sym^2g^.
\to\cdots$$
where $^*J^{(i)}_m$ denotes the sum of the
terms in $^*J^._m$ of {\it{tensor}} degree $i$
(i.e. products of $i$ factors) and
$(^*J^{(i)}_m)_{[-i}$ is its truncation in (total) degrees
$\geq -i$.
Now the duality pairing
$$g^.\times ^*g^.\to \C$$
extends to 'contraction' maps (analogous to interior
multiplication)
$$(\Gamma^*J^{(r)}_m)_{[-r}\otimes\sym^2\g^.\to
(\Gamma^*J^{r-1}_{m-1})_{[-r+1}
\otimes g^..$$
Thanks to the alternating nature of the bracket on $g^.$,
it is easy to check
that these maps together yield a map of double complexes
$$ \Lt_m(\sym^2g^.)\to \Lt_{m-1}(g^.)[-1]. $$

Now recall the map
$$\sym^2\Lt_m(g^.)\to \Lt_m(\sym^2g^.)$$
as in (5.4).  Composing, we get a map of double complexes
$$b_m:\sym^2\Lt_m(g^.)\to \Lt_{m-1}(g^.)[-1], $$
which induces a map on the respective
truncations, whence a map on cohomology
$$H^2((\sym^2\Lt_m(\g))_+)\to H^1(\Lt_{m-1}(\g)_+)=
H^1(L_{m-1}(\g)). $$
Note that $(\sym^2\Lt_m(\g))_+=\sym^2(\Lt_m(\g)_+)=
\sym^2L_m(\g)$ because these are complexes in nonpositive
vertical degrees.
Then define the bracket  pairing
$$B_m:\bigwedge^2\Theta_m(\g)\to \Theta_{m-1}(\g)$$
as the induced map
$$\bigwedge^2 H^1(L_m(\g))\to H^2(\sym^2(L_m(\g)))\to
H^1(L_{m-1}(\g)).$$
\proclaim{Theorem 6.1} \item{(i)} The above pairings (6.9)
yield a compatible
sequence of natural homomorphisms
$$A_m:\Theta_m(\g)\to \Der (R_{m+1}(\g),R_{m}(\g)).$$
$A_1$ is always an isomorphism.
\item{(ii)} Via $A_m$, the commutator of derivations is given by
$$[A_{m+1}(u),A_{m+1}(v)]=A_m(B_m(u\wedge v)), \forall
u,v\in\Theta_{m+1}(\g).$$
\item{(iii)}
The induced pairing $\hat{B}
=\lim\limits_{\leftarrow}B_m$ on $\hat{\Theta}(\g)
=\lim\limits_{\leftarrow}
\Theta_m(\g)$ turns it into a Lie algebra.
If $\g$ has unobstructed deformations, then the induced map
$$\hat{A}=\lim\limits_{\leftarrow}A_m:\hat{\Theta}(\g)
\to\Der (\Hat{R}(\g))$$
is a Lie isomorphism.
\endproclaim
\demo{proof}(i) We first show $\Theta_m$ acts on $R_{m+1}$ (we drop
the $\g$ for convenience) as derivations, i.e. that
$$<A_m(u),fg>=g<A_m(u),f>+f<A_m(u),g>,
\ \ \forall f,g\in R_{m+1}, u\in\Theta_m.$$
 This results from the commutative diagram
$$\matrix L_m(\g)\times L_{m+1}(\C)\times L_{m+1}(\C)&
\longrightarrow
&L_m(\C)\times L_{m}(\C)[-1]\\
\id\times\mu_{m+1}\downarrow&&\mu_m\downarrow\\
L_m(\g)\times L_{m+1}(\C)&\overset <.>\to
\longrightarrow &L_m(\C)[-1]\endmatrix\tag 6.10 $$
where $\mu_m:L_m(\C)\times L_m(\C)\to L_m(\C)$
is the multiplication mapping
as in (5.4), which yields the multiplication in $R_m$, and which is
simply induced by (graded) exterior multiplication
in $\bigwedge\limits^.(^*g)$, $<.>$ is the
pairing (6.9), and the top horizontal arrow in induced by $<.>$ via
the derivation rule, i.e. $u\times f\times g\mapsto
<u,f>\times g+<u,g>\times f$.
Commutativity of this diagram is immediate from the definitions.\par
Next we show $A_m$ is $R_m$-linear.
This again follows from the (easily checked)
commutativity of a suitable diagram, namely
$$\matrix L_m(\C)\times L_m(\g)\times
L_{m+1}(\C)&\overset \id\times <.>\to\longrightarrow&
L_m(\C)\times L_m(\C)[-1]\\
\mu'\times\id\downarrow&&\downarrow\mu_m\\
L_m(\g)\times L_{m+1}(\C)&\overset <.>\to
\longrightarrow &L_m(\C)[-1]\endmatrix\tag 6.11$$
where $\mu'$ is the multiplication mapping
$L_m(\C)\times L_m(\g)\to L_m(\g)$ which
induces the $R_m$-module structure on $\Theta_m$.\par
Note that $A_1$ is just a map
$$H^1(\g)\to\m_{R_1(\g)}^*=H^1(\g),$$
and it is immediate from the definitions
that this is just the identity.\par
(ii) To be precise, what is being asserted here is that for all
$u,v\in\Theta_{m+1},$ if $u',v'$ are the induced elements
in $\Theta_m$, then
$$A_{m}(u')\circ A_{m+1}(v)-A_{m}(v')\circ
A_{m+1}(u)=A_m(B_m(u\wedge v)).$$
This, in turn, results from the
commutativity of the following diagram
$$\matrix \sym^2L_m(\g)\times L_{m+1}(\C)&
\overset b_m\times\bar{\id}\to\longrightarrow&
L_{m-1}(\g)\times L_m(\C)[-1]\\
\id\times<.>\downarrow&&\downarrow <.>\\
L_{m-1}(\g)\times L_m(\C)[-1]&\overset <.>\to
\longrightarrow& L_{m-1}(\C)[-2]\endmatrix$$
where the left vertical arrow is induced by $<.>$ again via
the derivation rule, i.e. $uv\times w
\mapsto u\times <v.w>+v\times <u.w>$.
\par
(iii) The fact that $\hat{\Theta}(\g)$
is a Lie algebra amounts to the Jacobi identity. To
verify this, note that $b_m$ induces via the derivation rule a map
$$\sym^3L_m(\g)\to L_m(\g^.)\otimes
L_{m-1}(\g)[-1]\to \sym^2L_{m-1}(\g)[-1].$$
Then the Jacobi identity amounts to
the vanishing of the composite of this map and
$$b_{m-1}:\sym^2L_{m-1}(\g)[-1]\to L_{m-2}(\g)[-2].$$
This may be verified easily.\par
Finally in the unobstructed case, clearly both
$\hat{\Theta}(\g)$ and
$\Der (\Hat{R}(\g))$ are free $\Hat{R}(\g)$-modules,
and since $A_1$ is an isomorphism
it follows that so is $\hat{A}$.\qed
\enddemo
\remark{Elaboration 6.2} In term of cocycles, we may
describe $\Theta_1(\g)$ as follows.
Set $V=H^1(\g)$ which we view as a subspace of $\Gamma (\g^1)$.Then
$$\Theta_1(\g)=\{ (a,\sum b_i\otimes c^*_i)\in
V\oplus\g^1\otimes V^*|
^tbr (a)=\sum\delta (b_i)\otimes c_i^*\}$$
where $^tbr$ is the adjoint of the bracket, defined by
$$^tbr (a)=\sum b_i\otimes c_i^*\,$$
$$[a,x]=\sum <c_i^*.x> b_i\ \ \forall x\in V.$$
Thus  the condition defining $\Theta_1(\g)$ is
$$[a,x]=\sum <c_i^*.x>\delta (b_i)\ \ \forall x\in V.$$
Now the bracket
$$[.,.]:\bigwedge\limits^2\Theta_1(\g)\to \Theta_0(\g)=V$$
is given by
$$[(a,\sum b_i\otimes c_i),(a',\sum b'_i\otimes c^{'*}_i)]=
\sum <c^{'*}_i.a>b'_i-\sum <c^{*}_i.a'>b_i.$$
Note that neither sum is $\delta -$ closed,
but the difference is  because
$$\sum <c^{'*}_i.a>\delta (b'_i)-\sum <c^{*}_i.a'>b_i
=[a',a]-[a,a']=0$$
(recall that the bracket is symmetric on $\g^1$).
\endremark
\heading 7. Differential operators\endheading
We shall require an extension of the results the the last
section from the case
of derivations on $R_m(\g)$ itself to that of differential
operators on 'modular' $R_m(\g)$-modules (those that come from
$\g-$ modules).To this end we will construct,
given a dgla $g^.$ and $g^.-$modules $G_1,G_2$,
 complexes $LD^m_k(G_1,G_2)$
whose cohomology will act as $m$th order differential operators
from $H^.(M_{m+k}(G_2)$ to $H^.(M_k(G_1))$ and will coincide
with the module of all such operators, i.e.
$D^m(H^.(M_{m+k}(G_2),H^.(M_k(G_1))$, under favorable circumstances
('no obstructions').This will apply in particular to
an admissible Lie pair $(\g, E)$
on $X$ with suitable (dgla, dg-module) resolution
$(\g^.,E^.)$, by taking as usual
$$g^.=\Gamma (\g^.),$$

$$G^.=\Gamma(E^.).$$\par
To begin with, set, for any $g^.$-modules $G_1,G_2$,
$$K_m(G_1,G_2)=K_m(g^.,G_{1 triv}\otimes ^*G_2),$$
where $G_{1 triv}\otimes ^*G_2$ is $G_1\otimes ^*G_2$
as a complex but with $g^.$ acting  through the $^* G_2$
factor only. Note that as a complex, we may identify
$K_m(G_1,G_2)=G_1\otimes K(g^.,^*G_2)$ .
A fundamental observation is the following
\proclaim{Lemma 7.1} Let $G_1, G_2$ be $g^.-$modules.
Then the duality pairings between
$g^.$ and  $^*g^.$ and $G_2$ and $^*G_2$
extends to a pairing
$$K_m(G_1, G_2)\times L_m(G_2)\to G_1.$$\endproclaim
\demo{proof} There is clearly no loss of generality in
assuming $G_1=\C$ with trivial $g-$action.
Write these complexes schematically as
$$K:=K_m(g^.,^*G_2)\ \ \cdots\bigwedge\limits^2g\otimes ^*G_2
\to g\otimes ^*G_2\to ^*G_2,$$
$$L:=L_m(g^.,G_2)\ \ G_2\to^*g\otimes G_2\to
\bigwedge\limits^2(^*g)\otimes G_2\cdots.$$
Then we have
$$(K\otimes L)_0=\bigoplus\limits_{i=0}^m\bigwedge\limits^ig\otimes
\bigwedge\limits^i(^*g)\otimes ^*G_2\otimes G_2.$$
We map this to $\C$ in the obvious way by contracting
together all the $g$ and $^*g$ factors and likewise for
$^*G_2$ and $G_2$.
What has to be proved is that
this yields a map of complexes $K\otimes L\to \C,$
i.e. that the composite
$$(K\otimes L)_{-1}\overset{\delta_{K\otimes L}}\to\to
(K\otimes L)_0\to \C$$
vanishes, in other words that for each $i=1,...,m$
the composite
$$\bigwedge\limits^ig\otimes\bigwedge\limits^{i-1}(^*g)\otimes ^*G_2
\otimes G_2
\overset{\delta_K\otimes\id\oplus \id\otimes\delta_L}\to\to
\bigwedge\limits^{i-1}g\otimes\bigwedge\limits^{i-1}(^*g)\otimes
^*G_2\otimes G_2
\bigoplus
\bigwedge\limits^ig\otimes\bigwedge\limits^{i}(^*g)\otimes
^*G_2\otimes G_2$$$$
\to ^*G_2\otimes G_2\to \C$$
is zero. Now it is easy to see from the definitions
(compare the proof of Lemma 5.8) that it
suffices to check this for $i=1.$ So pick an element
$$v\times a^*\times a\in g\times ^*G_2\times G_2.$$
Its image under the first map has two components,
the first of which is $<v,a^*>\times a$ where $<\cdot,\cdot>$
denotes the action, while the second component has the form
$$v\times a^*\times\sum w^*_j\otimes b_j $$
where the sum denotes the cobracket of $a$, defined by
$$\sum <w^*_j.y>b_j=<y,a>,\forall y\in g,$$
where $<\cdot.\cdot>$ denotes the duality pairing. Clearly the image
of this second component in $\C$ (i.e. its trace) is just
$<a^*.<v,a>>$. However by definition of the dual action
we have
$$<a^*.<v,a>>=-<<v,a^*>.a>.$$
Thus the image of $v\times a^*\times a$ in $\C$ is zero, as claimed.
\qed
\enddemo

Next, recall by Lemma 5.8 that $K_m(^*G_2)$ is a $g^.$-module,
hence so is $K_m(G_1,G_2)=G_1\otimes K_m(^*G_2).$
This gives rise to a complex
$$L_k(g^., K_m(G_1, G_2))=:LD_k^m(g^.,G_1,G_2).$$
When $g^.$ is understood, we may denote the latter by
$LD_k^m(G_1,G_2),$ and when $G_1=G_2=G$ the same may also be denoted
by $LD_k^m(G)$
 From (5.4) and Lemma 7.1
 we deduce a pairing
 $$LD^m_k(g^.,G_1, G_2))\times L_m(g^.,L_k(g^.,G_2))\to
L_k(G_1),$$ hence by Lemma 5.8(ii) we get a pairing
$$LD^m_k(g^.,G_1, G_2))\times L_{m+k}(g^.,G_2)\to
L_k(G_1)\tag 7.1$$

Our next goal is to show that,
via this pairing, we may, at least under
favorable circumstances, identify $H^{j-i}(LD_k^m(g^.,G_1,G_2)$
with the $k$-jet of the $m-$th order differential operators
on the $R_{m+k}(g^.)$-module corresponding to $H^{i}(G_2)$
 with values in $H^{j}(G_1)$, provided these are
 the unique
nonvanishing respective cohomology groups.
In fact, it will be convenient to prove the stronger result
saying that this assertion
essentially holds already 'on the level of complexes'.
To explain what that means, note that via the pairing (5.4),
$L_m(\C)$- and likewise $L_m(A)$ for any $\C$-algebra $A$-
forms a 'ring complex', i.e. a ring object in the category
of complexes;
this ring structure is the one that induces the ring structure on
$R_m(g)=\H^0(\lmc)$.
 Moreover, for any $g$-module $G$, $L_m(G)$
is an $L_m(\C)-$module. There is an evident notion of
$\lmc-$linear map of $\lmc$-modules, and any $g-$linear map
$G_1\to G_2$ induces such a map $L_m(G_1)\to L_m(G_2).$
Likewise, the natural pairing
$$L_m(G)\times L_m(^*G)\to\lmc$$
is $\lmc-$linear.
\par
Given this, the notion of differential operators of any order
(over $\lmc$) can be defined inductively: given complexes
$D,M,N$ of $\lmc-$modules and a pairing
$$a: D\times M\to N,$$
$a$ is said to be of differential order $\leq m$ in the
$M$ factor if the composite map
$$\lmc\times D\times M\to N,$$
$$(v,d,m)\mapsto a((vd),m)-a(d,(vm))$$
is of differential order $\leq m-1$ in the $M$ factor.
\proclaim{Lemma 7.2} The pairing (7.1) is of differential order
$\leq m$ in $L_{m+k}(G_2).$\endproclaim
\demo{proof} By induction on $m$, of course. For $m=0$ the result
is clear (and was already noted above). For the induction
step, it suffices to show that the map
$$LD^m_k(G_1,G_2)\times L_{m+k}(\C)\to LD^m_k(G_1,G_2)$$
given by (premultiplication)-(postmultiplication) factors
through $LD^{m-1}_k(G_1,G_2)$. As for the premultiplication map,
it is induced by the $\lmc-$module structure on $K_m(^*G_2),$
i.e. the natural map (cf. Lemma 5.8(i))
$$\lmc\times K_m(^*G_2)\to K_m(^*G_2).$$
Tensoring by $G_1,$ applying $L_k$ and using Lemma 5.8(ii) we get a map
$$L_{m+k}(\C)\times LD^m_k(G_1,G_2)\to LD^m_k(G_1,G_2)$$
that is the premultiplication map. This map clearly factors
through $L_{m}(\C)\times LD^m_k(G_1,G_2)$. It
is essentially
obtained by contracting together some $g$ and $^*g$ factors
and exterior-multiplying others;
in particular the induced map on any term involving
$\bigwedge\limits^m(g)$ going to a similar term cannot
involve any contraction, hence is simply given by exterior-
multiplying the factor from $\lmc$ by the one from
$LD^m_k(G_1,G_2)$. It is easy to see that similar comments
apply to the postmultiplication map. Thus the two induced
map (from pre and post)
between terms involving $\bigwedge\limits^m(g)$
agree, and consequently the difference (pre)-(post)
goes into $LD^{m-1}_k(G_1,G_2)$, which proves the Lemma.\qed\enddemo
\remark{Remark 7.2.1} As was observed in the course of the
proof, $LD^m_k(G_1,G_2)$ has the structure of $L_{m+k}(\C)-$bimodule,
corresponding to the pre-post-multiplications. This is analogous,
and closely related to the bimodule structure on the space of
differential operators $D^m(M_1,M_2)$ between a pair
of modules.\endremark
\par
Next we will construct a pairing that will yield
the {\it composition} of differential operators.
\proclaim{Lemma 7.3} For any $g-$modules $G_1,G_2,G3$ and
natural numbers $m,k,j,n$ with
 $k\geq j-m\geq 0$,there is  a natural pairing of $g-$
 modules
$$LD_k^m(G_1,G_2)\times LD_j^n(G_2,G_3)\to LD^{m+n}_{j-m}(G_1,G_3).
\tag 7.2$$
\endproclaim\demo{proof} There is clearly
no loss of generality in assuming $k=j-m$. Then using
Lemma 5.8 we are easily reduced to the case $j=m$ where
the point is to construct a $g-$linear pairing
$$G_1\otimes K_m(^*G_2)\times L_m(G_2\otimes K_n(^*G_3))
\to G_1\otimes K_{m+n}(^*G_3).\tag 7.3$$
There is obviously no loss of generality in assuming
$G_1=\C.$ Then the LHS is a direct sum of terms of the form
$$\bigwedge\limits^ig\otimes ^*G_2\times
\bigwedge\limits^j(^*g)\otimes\bigwedge^kg\otimes G_2\otimes G_3$$
which has degree $i+k-j$. We map this term to zero
if $i+j-k<0$, and otherwise to
$$\bigwedge\limits^{i-j+k}g\otimes ^*G_3=K_{m+n}^{i-j+k}(^*G_3)$$
in the standard way, by contracting away all the
$^*g$ factors against the $g's,$
as well as $^*G_2$ against $G_2$. If we can
prove this is a map of complexes then $g-$linearity
comes for free, due to the $g-$linearity of contraction.\par
Now to prove we have a map of complexes one may reduce as in
the proof of Lemma 5.8 to the case $i=k=1,j=0$ and commutativity
of the following diagram
$$\matrix g\otimes ^*G_2\otimes g\otimes G_2\otimes ^*G_3&\to&
[g\otimes ^*G_2\times G_2\otimes ^*G_3]^{\oplus 2}\oplus
g\otimes ^*G_2\otimes ^*g\otimes g\otimes G_2\otimes ^*G_3\\
\downarrow&&\downarrow\\
\bigwedge\limits^2g\otimes ^*G_3&\to&g\otimes ^*G_3\endmatrix\tag 7.4$$
where the top map is of the form\par
{($g$-action on $^*G_2$, $g$-action on $^*G_3$, $g$-coaction on
$g\otimes G_2\otimes^*G_3$)}.\par
Given an element
$$v_1\times a^*\times v_2\times a\times b^*\in
g\otimes ^*G_2\otimes g\otimes G_2\otimes ^*G_3,$$
its image going counterclockwise is clearly
$$<(v_1\wedge v_2),<a.a^*>b^*>=$$
$$<a.a^*>(v_1\times <v_2,b^*>
-v_2\times <v_1,b^*>-[v_1,v_2]\times b^*).\tag 7.5$$
On the other hand, the image of this element under the top
map is
$$(<v_1,a^*>\times v_2\times a\times b^*, v_1\times a^*\times
<v_2,b^*>\times a, v_1\times a^*\times\cob(v_2\times a\times b^*)).$$
Now from the definition of $\cob$, the fact that it acts as a derivation,
plus the definition of the dual action, it is elementary to
verify that the image of the latter element under the right
vertical map coincides with (7.5), so the diagram commutes.\qed

\enddemo
\proclaim{Lemma 7.4} Via the action pairing (7.1), the 'composition'
pairing (7.2) corresponds to composition of operators.\endproclaim
\demo{proof} Our assertion
means that
$$<<d_1,d_2>,a>=<d_1,<d_2,a>>,$$
$$\forall d_1\in
LD^m_k(G_1,G_2), d_2\in LD^n_j(G_2,G_3), a\in L_r(G_3),$$
assuming $r\geq j-m\geq 0$ (and abusing $<\ >$ to denote the
various pairings involved),
which amounts to commutativity of
the obvious diagram
$$\matrix LD^m_k(G_1,G_2)\times LD^n_j(G_2,G_3)
\times L_r(G_3)&\to&LD^{m+n}_{j-m}(G_1,G_3)\times L_r(G_3)\\
\downarrow&&\downarrow\\
LD^m_k(G_1,G_2)\times L_j(G_2)&\to&L_{j-m}(G_1).\endmatrix\tag 7.6$$
 Now all the maps involved
are essentially given by exterior multiplication and
contraction, so commutativity of (7.6) follows from the
associativity of exterior multiplication.\qed

\enddemo
 In
particular, taking $G_1=G_2=G_3=G$
we get a (composition) pairing, whenever $k\geq m,$
$$LD_k^m(G)\times LD_k^n(G)\to LD_{k-m}^{m+n}(G).$$
It is easy to see by sign considerations
as in the proof of Lemma 7.2 that the
'commutator' associated to this pairing takes values in
$LD_{k-\max(m,n)}^{m+n-1}(G).$ In particular, we get a skew-symmetric
pairing
$$B_k:\bigwedge\limits^2LD^1_k(G)\to LD^1_{k-1}(G).$$
\proclaim{Lemma 7.5}
Under ${B_\infty}=\lim\limits_\leftarrow B_k$,
${LD^1_\infty(G)}=\lim\limits_\leftarrow LD^1_k(G)$ is a Lie algebra
object in the category of complexes,
and admits a natural representation on
$L_\infty(G)=\lim\limits_\leftarrow L_k(G)$.\endproclaim
\demo{proof} Most of this has been proved above.
The only remaining point is the Jacobi identity
for $B_\infty$, which can be proved as in the case of the
trivial module $G=\C$ (cf. Theorem 6.1).\qed\enddemo
The pairings discussed above naturally induce analogous
pairings on cohomology groups. This leads to the following
Theorem 7.6 . First some
notation and terminology. For any $g-$module $G, k\leq\infty ,$ set
$$H^i(G,k)=H^i(L_k(G))$$
As we have seen if $(g,G)$ comes from sheaves $(\g, E)$
on $X$ then
this coincides with the sheaf cohomology
$H^i(X, M_k(\g,E))$, i.e. the $k-$universal
$g-$deformation of $H^i(X, E).$ We will say that $G$ is
$strongly\ i- unobstructed$ if for all $v\in g^1$
that is $ \delta-$closed
(i.e. $\delta(v)=0$) and all $a\in G^i$ (closed or not),
we have that $<v,a>$ is exact; we will say that $g$ itself
is {\it{strongly unobstructed}} if it is strongly 1-unobstructed
in the adjoint representation. It is easy to see that if
$g$ is strongly unobstructed then $R_\infty(g)$ is regular
(i.e. smooth) and that if $G$ is strongly $i$-unobstructed then
$H^i(G,\infty)$ is $R_\infty(g)$ -free. Also, it is obvious that
if $G$ is $i$-equicyclic then it is strongly $i$-unobstructed.

\proclaim{Theorem 7.6} Let $G_1,G_2,G_3$ be modules over
the dgla $g$ with $H^{\leq 0}(g)=0.$ Then
\item{(i)} there is a natural pairing, for any $0\leq k\leq n-m$
$$H^{j-i}(LD^m_k(G_1,G_2))\times H^i(G_2,n)\to H^j(G_1,k)$$
which induces a map
$$A:H^{j-i}(LD^m_k(G_1,G_2))\to D^m_{R_n(g)}
( H^i(G_2,n), H^j(G_1,k);$$
\item{(ii)} there is a natural pairing, for any $0\leq j-m\leq k$,
$$C:H^i(LD^m_k(G_1,G_2))\times H^j(LD^n_j(G_2,G_3))\to
H^{i+j}(LD^{m+n}_{j-m}(G_1,G_3),$$
via which $A$ corresponds to composition of operators;
in particular there are  natural Lie algebra
structures on $H^0(LD^1_\infty(G))$ and $H^0(LD^\infty_\infty(G))$
with  representations on $H^i(G,\infty)$ for all $i$;
\item{(iii)} if $g$ is strongly unobstructed and $G_1$ and $G_2$
are equicyclic of degrees $i,j$ respectively, then the map
$$A_\infty:H^{i-j}(LD^m_\infty(G_1,G_2))\to D^m_{R_\infty(g)}
(H^j(G_2,\infty),H^i(G_1,\infty))\tag 7.7$$
is an isomorphism for all $m$.
\endproclaim
\demo{proof} Items (i) and (ii) follow directly from the
results above. We prove (iii). Clearly the target of
$A_\infty$ , with respect to its left
(postmultiplication) structure, is a free
module with fibre
$$H^i(G_1)\otimes D^m_0(H^j(G_2),\C)=H^i(G_1)\otimes H^j(G_2,m)^*.$$

As for the source, note that
$K_m(^*G_2)$ is strongly $(-j)$-unobstructed and has no cohomology
in degree $<-j.$ Consequently, $G_1\otimes K_m(^*G_2)$
is strongly $(i-j)$-unobstructed and
$$H^{i-j}(G_1\otimes K_m(^*G_2))=H^i(G_1)\otimes H^{-j}(K_m(^*G_2))
=H^i(G_1)\otimes H^j(G_2,m)^*.$$
By definition, the latter is precisely the fibre of
$H^{i-j}(LD^m_\infty(G_1,G_2))$ with respect to its postmultiplication
module structure (which structure we now know is free,
thanks to unobstructedness). Thus
the source and target of $A_\infty$ have isomorphic fibres; moreover
it is easy to see, for instance by considering the other
(right or premultiplication) structure that $A_\infty$ induces
an isomorphism. But clearly a linear map of free modules
over a local ring inducing an isomorphism on fibres is itself an
isomorphism, proving our assertion.\qed
\enddemo
\proclaim{Corollary 7.7} If $G$ is an $i-$equicyclic module and
$g$ is strongly unobstructed then the Lie algebra $H^0(LD^1_\infty
(G))$ is canonically isomorphic to $D^1_{R_\infty(G)}(H^i(G,\infty))$
\qed\endproclaim
In particular, in the geometric situation with
$(\g,E)$ an admissible pair, $\g$ unobstructed and $E$ i-equicyclic,
we get a canonical recipe for the
Lie algebra which is the
formal completion of $D^1_\M(\SH)$ where $\SH$ is the sheaf on
the moduli space $\M$ associated to the unique nonvanishing
cohomology group $H^i(E)$.\par
\remark{Elaboration 7.8}
Let us write down the bracket
pairing $B_1$ in terms of cocycles.
 This comes about by considering the diagram
$$\matrix g\otimes\ G\otimes ^*G&\overset b\to\longrightarrow&
G\otimes ^*G\\
^tb\downarrow&&^tb\downarrow\\
^*g\otimes g\otimes\ G\otimes ^*G&\overset b\to\longrightarrow&
^*g\otimes\ G\otimes ^*G\endmatrix\tag 7.8$$
where the maps $b$ are induced by the action of $g$ on $^*G$, while
the maps $^tb$ are induced by the transpose of the $g$ action on
$G$. A cocycle for $LD^1_1(G)$ is a 4-tuple
$(\phi,\psi,\phi',\psi')$of cochains of the four complexes in (7.8)
such that
$$\p (\phi)=0$$
$$b(\phi)=\p (\psi)$$
$$^tb(\phi)=\p (\phi')$$
$$b(\phi')+^tb(\psi)=\p (\psi').$$
The pairing
$$B_1:\bigwedge\limits^2(H^0(LD^1_1(G)))\to H^0(LD^1_0(G)$$
is given by
$$B_1((\phi_0,\psi_0,\phi'_0,\psi'_0)\wedge
(\phi_1,\psi_1,\phi'_1,\psi'_1))=(\phi_2,\psi_2)$$
where
$$\psi_2= [\psi_0,\psi_1]+<\phi_0,\psi'_1>-<\phi_1,\psi'_0>$$
$$\phi_2=<\phi_0,\phi'_1>-<\phi_1,\phi'_0>$$
(compare Elaboration 6.2).
Here $[\ ]$ is the usual commutator on $G\otimes ^*G$
while $<\ >$ is the pairing induced by $[\ ]$ and the duality
between $g$ and $^*g$.
To check that this is indeed a cocycle, we compute:
$$\p (\psi_2)=[\p (\psi_0),\psi_1]-[\psi_0,\p (\psi_1)]
-<\phi_0,\p (\psi'_1)>+<\phi_1,\p (\psi'_0)>$$
$$=[b(\phi_0),\psi_1]-[\psi_0,b(\phi_1)]
-<\phi_0,b(\phi'_1)+^tb(\psi_1)>+
<\phi_1,b(\phi'_0)+^tb(\psi_0)>$$
$$=[b(\phi_0),\psi_1]-[\psi_0,b(\phi_1)]
-<\phi_0,b(\phi'_1)>+<\phi_1,b(\phi'_0)>-
[b(\phi_0),\psi_1]+[\psi_0,b(\phi_1)]$$
$$=-<\phi_0,b(\phi'_1)>+<\phi_1,b(\phi'_0)>$$
$$=b(<\phi_0,\phi'_1>-<\phi_1,\phi'_0>)=b(\phi_2).$$
\endremark
Analogous formulae may be given for the bracket 'action' of
$LD^1_k(G)$ on $LD^m_k(G)$. These actions being compatible
for different $m$, there is an induced action
on $LD^m_k(G)/L^{m-1}_k(G)=L_k(\bigwedge\limits^m g\otimes
G\otimes ^*G)[m]$.
In particular, we get a pairing
$$LD^1_1(G)\times L_1(\bigwedge\limits^2
g\otimes G\otimes ^*G)[2]\to
(\bigwedge\limits^2 g\otimes G\otimes ^*G)[2]$$
Now note the natural map
$$L_1(\bigwedge\limits^2 g\otimes G\otimes ^*G)[1]\to LD^1_1(G)$$
which is induced by the map $\bigwedge\limits^2
g\otimes G\otimes ^*G[1]\to
K_1(g,G\otimes ^*G)$ that is part of the complex
$K_2(g,G\otimes ^*G)$.
Hence we get a pairing
$$L_1(\bigwedge\limits^2 g\otimes G\otimes ^*G)[1]\times
L_1(\bigwedge\limits^2 g\otimes G\otimes ^*G)[2]\to
(\bigwedge\limits^2 g\otimes G\otimes ^*G)[2]$$
i.e. a (symmetric) bracket pairing
$$\Sym^2(L_1(\bigwedge\limits^2 g\otimes G\otimes ^*G)[1])
\to \bigwedge\limits^2 g\otimes G\otimes ^*G[1].$$
This pairing has the following interpretation.
Suppose $\M$ is a locally fine moduli space with
Lie algebra $\tilde{\g}$ on $X\times \M$ as above
and $\SH$ is the locally free $\O_\M-$sheaf $R^ip_{\M*}(\SH)$
for a suitable $\g-$module $E$ on $X\times\M$ (assuming this is
the only nonvanishing derived image).
Then as in Example 1.1.2 C we get
a heat atom
$$(\D^1_\M(\SH),\D^2_\M(\SH))$$
on $\M$, hence a Lie bracket on the (shifted) quotient
$\sym^2(T_\M)\otimes \SH^*\otimes \SH[-1]$. This bracket
can be identified 'fibrewise' with the map induced
by the pairing (7.1).

 \vfill\eject
 \heading 8. Connection Algebra
 \endheading Our
purpose in this section is to construct, for a given
representation $(\g,E)$, a canonical 'thickening' $\k(\g,E)$ of
$\g$ which is another Lie algbera which
 acts on $E$, such that the sections of $E$ extend
canonically over the universal deformation associated to
$\k(\g,E)$.\par Our construction refines and generalizes one
in first-order deformation theory due
to Welters [W] and Hitchin [Hit, Thm 1.20].
They noted that given a line
bundle $L$ on a compact complex manifold $X$, together with a
holomorphic section $s\in H^0(L)$,
1-parameter deformations of the
triple $(X,L,s)$ are parametrized by $\H^1$ of the complex
$$\D^1(L)\overset {\cdot s} \to \longrightarrow L.$$
Consequently, a family, in a suitable sense,
of such $\H^1$ cohomology
classes yields a {\it{connection}} $\nabla$ on the
'bundle of $H^0(L)$'s
(more precisely, it yields the covariant
derivative $\nabla\cdot s$).\par
Our construction, amongst other things,
extends that of Welters-Hitchin
from first-order to arbitrary $m$-th order
deformations. Applied in their
original context with $m$ at least 2,
it shows that the connection
$\nabla$ is automatically ${\underline{flat}}$,
a fact which could not be seen by
first-order considerations alone.\par
Now let $(\g,E)$ be an admissible pair,
with soft resolution $(\g^.,E^.,\p)$.
Then $\Gamma^*(E^.)\otimes E^.$ is a complex
(via tensor product of complexes)
and a $\g^.$-module (acting on the
$E^.$ factor only), which makes it
a \dg $\g^.-$module. There is a tautological map
$$\g^.\overset\delta\to\longrightarrow
\Gamma^*(E^.)\otimes E^.\tag 8.1$$
which is easily seen to be a derivation. Thus, (the mapping
cone of) (8.1) yields a \dgla, which we denote
$\k(\g,E)$. Note that $\k(\g,E)$ is
itself a \dg $\g^.-$module, and that we
have a natural dgla homomorphism
$$\k(\g,E)\to \g^.$$
Note also that if $H^{\leq 0}(\g)=0$, then we have
$$H^{\leq 0}(\k(\g,E))=0$$
if and only if $E$, that is,
$\Gamma E^.$, is $i$-equicyclic
for some $i$, in which case
we have an exact sequence
$$0\to H^i(E)\otimes
H^i(E)^*\to H^1(\k(\g,E))\to H^1(\g).$$\par
Similar constructions can be make purely algebraically.
Thus let $(g^.,G^.)$ be
a dg Lie representation.
We consider $^*G^.\otimes G^.$ as another dg
representation of $g^.$
(with differential as tensor product of complexes
and $g^.$-action on the
$G^.$ factor only), and note the graded derivation
$$g^.\overset\delta\to\longrightarrow\ ^*G^.\otimes G^..
\tag 8.2$$
Then (8.2) forms a dgla which we denote by $k(g^.,G^.)$,
and in which $g^i$ has degree $i$ and $^*G^i\otimes G^j$
has degree $i+j+1$. Thus
$$\Gamma\k(\g,E)=k(\Gamma\g^.,\Gamma E^.).$$
Obviously, $k(g^.,G^.)$ is a $g^.$-module; indeed the
canonical 'identity' element
$$I\in \ ^*G\otimes G$$ yields an inclusion
of $g^.-$modules
$$k(g^.,G^.)\subset LD^1_0(G)$$ (cf. \S 6).
Note that the $g^.$-action on $ ^*G^.\otimes G^.$
evidently extends to an
action of $k^.=k(g^.,G^.)$, by letting
$ ^*G^.\otimes G^.$ act trivially on itself.
Consequently we get for each $m\geq 1$ a complex
$L_m(k(g^.,G^.),\ ^*G^.\otimes G^.))$ which we write
schematically as a double complex
(with components which are themselves multiple complexes)
in the form

$$\matrix
&&&&\sym^2(^*G^.\otimes G^.)\otimes\ ^*G^.
\otimes G^.&\to&\ldots\\
&&&&\downarrow&\\
&&^*G^.\otimes G^.\otimes\ ^*G^.\otimes G^.&\to&
\ ^*g^.\otimes\ ^*G^.\otimes G^.
\otimes\ ^*G^.\otimes G^.&\to&\ldots\\
&&^*\delta\otimes\id\downarrow&&\downarrow&& \\
^*G^.\otimes G^.&\to&\ ^*g^.\otimes\ ^*G^.\otimes G^.&
\to&\bigwedge\limits^2
(^*g^.)\otimes\ ^*G^.\otimes G^.&\to&\ldots\endmatrix
\tag 8.3$$

Thus the $i-$th column in (8.3) is the complex
$\bigwedge\limits^i(^*k(g^.,G^.))\otimes\ ^*G^.\otimes G^.$.

\proclaim{Lemma 8.1} The identity element
$I\in\ ^*G^.\otimes G^.$ lifts canonically
to a compatible sequence of elements
$$I_m\in H^0(L_m(k(g^.,G^.),\ ^*G^.\otimes G^.)), m\geq 1.$$
\endproclaim
\demo{proof} Let $I_m$ be the cochain consisting
of the elements
$\sym^iI\otimes I$ in position $(i,i)$ in the above complex,
for all $0\leq i\leq m$. It is trivial
to check that this is a cocycle.\qed\enddemo
\proclaim {Theorem 8.2} In the situation of Theorem 5.1,
assume moreover that $E$
is equicyclic of degree $ i$. Then we have a
canonical isomorphism
(or 'trivialization')
$$H^i(M_m(\g,E))\otimes_{R_m(\g)}R_m(\k(\g,E))\simeq
H^i(E)\otimes_{\C}R_m(\k(\g,E))
\tag 8.4$$
Moreover, $R_m(\k(\g,E))$ is universal
with respect to this property,
i.e. given a deformation $E^\tau$ parametrized
by $S$ and an
$S$-isomorphism
$$H^i(E^\tau)\simeq H^i(E)\otimes S$$
lifting the identity on $H^i(E)$,
there is a canonical lifting of the Kodaira-Spencer
homomorphism of $\tau$ to a homomorphism
$R_m(\k(\g,E))\to S.$
\endproclaim
\demo{proof} Apply Lemma 8.1 to
$g^.=\Gamma (\g^.), G^.=\Gamma (E^.).$
Because $g^.$ acts trivially on $^*G^.$, we have
$$L_m(k(g^.,G^.),\ ^*G^.\otimes G^.))=
\Gamma^*(E^.)\otimes L_m(\Gamma (\k(\g,E)),\Gamma (E^.)).$$
 As $H^j(\Gamma^*(E^.))=0$ for $j\neq -i$, we have
$$H^i(L_m(k(g^.,G^.),\ ^*G^.\otimes G^.)))=
\hom (H^i(E), H^i(L_m(\Gamma (\k(\g,E)),\Gamma (E^.)))).$$
Clearly
$$H^i(L_m(\Gamma (\k(\g,E)),\Gamma (E^.))))
\simeq H^i(L_m(\g,E)))
\otimes_{R_m(\g)}R_m(\k(\g,E)),$$
and by Theorem 5.1 this is just $H^i(M_m(\g,E))
\otimes_{R_m(\g)}R_m(\k(\g,E))$,
so the element $I_m$ above yields the required trivialization (8.4).
\par
In terms of cocycles, this trivialization may be
seen as follows. Consider the universal $\k(\g,E)$-deformation over
$R=R_m(\k(\g,E))$. This may be represented by
$$\psi=(\phi, \sum t_j\otimes t_j^*)\in (\Gamma(\g^1)
\oplus\Gamma (E^i)\otimes\Gamma (E^i)^*)\otimes\m,
\ \m=m_R. $$
Letting $(s_k\in \Gamma (E^i))$ be a lift of
some basis of $H^i$ and $s_k^*$ be a lift of
the dual basis, we may write
the integrability condition $\p\psi=-\1/2[\psi,\psi]$ as
$$\p \phi =-\1/2 [\phi,\phi ],$$
$$\sum (\p t_j)\otimes t_j^* =
-\sum [\phi,t_j]\otimes t_j^*-\sum [\phi,s_k]\otimes s_k^*,
\tag 8.5$$
$$\sum t_j\otimes (\p t_j^*)=0.$$
Thus, we may assume that $\p t_j^*=0$
hence we may adjust notation so that
$t_j^*=s_j^*.$
Then we may write 8.5 in the form
$$\sum \p(s_j+t_j)\otimes s_j^*) +
\sum [\phi,s_j+t_j]\otimes s_j^*=0\tag 8.6$$
Recalling that the deformation $E^\phi$ of $E$
induced by $\phi$ is just $(E^.,\p+\ad\phi)$,
8.6 shows precisely that $\sum (s_j+t_j)\otimes s_j^*$
is a lift of $I=\sum s_j\otimes s_j^*$
to $E^\phi\otimes R$, yielding a canonical $R-$valued lift of
each $s_j$.\par
The latter description makes it easy to
establish the universality of $R(\k(\g,E))$, thus completing the
proof. Given $E^\tau/S$, a lifting of the identity
on $H^i(E)$ to an $S$-isomorphism
$H^i(E)\otimes S\simeq H^i(E^\tau) $ is given by
an element
$$\sum t_j\otimes s_j^* \in\Gamma(E^i)\otimes
\Gamma(E^i)^*\otimes\m_S$$
(i.e $s_j+t_j$ is a lifting of $s_j$),
and writing down the condition that $s_j+t_j$ is a cocycle for
$\partial+\ad \tau$ and computing as above shows precisely that
$$\rho=(\tau, \sum t_j\otimes s_j^*)$$
is an $S$-valued cocycle for $\k(\g,E)$, yielding the
desired homomorphism $R(\k(\g,E))\to S.$\qed\enddemo
For $m=1$ this result (or rather, its 'relative version' )
 generalizes
the Welters-Hitchin construction of connections (see [Hi],
Thm 1.20). For $m\geq 2$ the trivialization
we construct amounts to showing that this connection is flat.\par

\vfill\eject
\heading 9.  Relative deformations over a global base\endheading
Our purpose in this section is to discuss a more global
and relative generalization of the
notion of deformation, which occurs not just
over a (thickened) point (represented by an artin
local algebra), but over a global base, suitably
thickened. This is closely related -but not identical- to the notion
of {\it family}
or {\it variation} of deformations; the slightly
subtle difference is illustrated by the fact
that a 'family of trivial deformations' may well be nontrivial
as a relative deformation (such subtleties however occur only
in the presence of symmetries locally over the base and
globally along fibres).\par
To proceed with the basic definitions, let
$$f:X_B\to B$$
be a continuous mapping of Hausdorff spaces
with fibres $X_b=f^{-1}(b)$ and base
$B$ which we assume endowed with
a sheaf of local $\C-$algebras
$\O_B$.
A {\it{Lie pair}}
$(\g_B,E_B)$ on $X_B/S$ consists of a sheaf $\g_B$
of $f\inv\O_B$-Lie
algebras (i.e. with $f\inv\O_B-$ linear bracket),
a sheaf $E_B$ of $f\inv\O_B-$modules
with $f\inv\O_B-$linear
$g_B-$action. This pair is said to
be {\it{admissible}} if it
admits compatible soft resolutions $(\g_B^.,E_B^.)$
such that $\g_B^.$ is
a dgla and $E_B^.$ is a dg representation of
$\g_B^.$, and moreover,
 $\Gamma (\g_B^.),\Gamma(E^._B)$ may be
linearly topologized
so that coboundaries (and cocycles) are closed, and
the cohomology is of finite type as $\O_B$-module
(and in particular vanishes in almost all degrees).
Let's call such resolutions
$good$. Note that if $(\g_B^.,E_B^.)$ is an admissible pair
then for any $b\in B$ the 'fibre'
$$(\g_b,E_b):=(\g_B,E_B)\otimes (\O_{B,b}/\m_{B,b})$$
is an admissible pair on $X_b.$\par
Now let $\SS$ be an augmented  $\O_B-$algebra of finite type
as $\O_B-$module,
with maximal ideal $\m_\SS$ (below we shall also
consider the case where $\SS$ is an inverse limit of
such algebras, hence is complete noetherian rather
than finite type).
By a {\it relative} $\g_B-${\it {deformation of}}
$E_B,$ {\it{parametrized
by}} $ \SS$ we mean a
sheaf $E^\phi_B$ of $\SS$-modules on $X_B,$
together with a maximal
atlas of the following data
\item{-} An open covering $(U_\alpha)$ of $X_B$.
\item{-}$\SS$-isomorphisms
$$\Phi_\alpha:E^\phi|_{U_\alpha}\simto E|_{U_\alpha}
\otimes_{\O_B}\SS.$$
\item{-} For each $\alpha,\beta,$ a lifting of
$$\Phi_\beta\circ\Phi_\alpha^{-1}\in\
Aut (E|_{U_\alpha\cap U_\beta}\otimes_{\O_B}\SS)$$
to an element
$$\Psi_{\alpha,\beta}\in
\exp (\g_B\otimes\m_{\SS}(U_\alpha\cap U_\beta)).$$
If $\g_B$ acts faithfully on $E_B$ then the
$\Psi_{\alpha,\beta}$ are uniquely
determined by the $\Phi_\alpha$ and form a cocycle;
in general we require
additionally that the
$\Psi_{\alpha,\beta}$ form a cocycle.\par
Note that if $(\g_B,E_B)$
is admissible then, as in the absolute
case,  for any relative
deformation $\phi$ there is
a good resolution $(E^.,\p )$ of $E$
and a resolution of $E^\phi$
of the form
$$ E^0\otimes_{\O_B}\SS
\overset{\p+\phi}\to
\longrightarrow E^1\otimes_{\O_B}\SS\cdots\tag 9.1$$
with $\phi\in\Gamma(\g^1_B)\otimes\m_\SS.$
We call such a resolution a {\it{good resolution}}
of $E^\phi$.
\remark{Example 9.1}(i) Let $E$ be a
vector bundle on the complex
manifold $X=X_B=B$
and let $\g=\gl (E).$ Let
$$P^m=P_X^m\O_{X\times X}/I_\Delta^{m+1},$$
which is naturally an $\O_X-$algebra via
the first coordinate projection $p_1$.
Likewise the $m-$th jet bundle
$$P^m(E)=p_{1*}(p_2^*(E)\otimes P^m)$$ is a
$P^m-$module and hence a $\g$-deformation
of $E$ parametrized by $P^m$
over $X$. Locally over the base $B=X,$
this deformation is obviously
trivial, but it is in general nontrivial
as relative deformation.
To obtain a good resolution of $P^m(E)$,
note that
$E$ admits a $\bp -$connection
(e.g. a Hermitian connection),
whose curvature is of type (1,1), i.e. trivial
on the (1,0) tangent directions,
hence yields a $C^\infty$ isomorphism
$$P^m(E)\sim P^m\otimes E,$$
hence the Dolbeault resolution of $P^m(E)$ is a
 good resolution as in (9.1).\par
More generally, $P^m(E)$ has a
structure of $\g-$deformation
for any $\O_X$-locally
free Lie subalgebra
$$\g\subseteq\gl (E)$$ such that
$E$ admits a $\g-structure$ (or 'reduction
of the structure algebra to $\g$'). To recall
what that means, let
$$G(E)=ISO (\C^r,E), r=\rk (E)$$
be the associated principal bundle, i.e. the open subset of
the geometric vector bundle $\hom (\C^r,E)$
consisting of fibrewise isomorphisms, with the obvious action
of
$GL_r$. Let $\D (E)$ be the sheaf of
$GL_r$-invariant vector fields on $G(E),$
which may also be identified as the sheaf of
'relative derivations' of $(E,\O_X)$,
consisting of pairs $(v,a), v\in T_X,
a\in Hom_\C(E,E)$ such that
$$a(fe)=fa(e)+v(f)e,\ \forall f\in \O_X,e\in E.$$
Note that $\D(E)$ is an extension of Lie algebras
$$0\to\gl(E)\to\D(E)\to T_X\to 0\tag 9.2$$
Then a $\g-$ structure on $E$ is a Lie subalgebra
$\hat{\g}\subseteq\D(E)$ which
fits in an exact sequence
$$\matrix 0&\to&\g&\to&\hat{\g}&\to&T_X&\to &0\\
           &   &\cap&&\cap&&\|&&\\
          0&\to&\gl(E)&\to&\D(E)&\to&T_X&\to&0.\endmatrix$$
Note that in this case a maximal
integral submanifold $\hat{G}$ of
$\hat{\g}$ yields a principal
subbundle of $G(E)$
with structure group $G=\exp (\g)$ and conversely
such a principal subbundle with Lie algebra $\g$
yields a $\g-$structure.
Clearly a $\g-$structure on $E$ yields a structure
of $\g-$ deformation
on $P^m(E)$ parametrized by $P^m$, for any $m$,
and as above this admits a good
(Dolbeault) resolution. We denote this deformation by
$P^m(E,\g)$.\par
Similarly, if $f:X_B\to B$ is any smooth morphism of complex
manifolds, and $E_B$ is a vector bundle on $X_B$ with a
relative $\g_B$-structure, then there is a relative $\g_B$-
deformation parametrized by $P^m_B$. We denote this deformation
by $P^m(E_B,\g_B)/B$ or by $P^m(E_B)/B$ if $\g_B$ is understood.
Intuitively, it represents the family of $m-$th order
deformations
$$E_B|_{f\inv (N_m(b))}=
E_B\otimes (\O_B/\m_{b,B}^{m+1}),\ b\in B,$$
where $N_m(b)=\Spec (\O_B/\m_{B,b}^{m+1})$ is the $m-$th order
neighborhood of $b$ in $B$.
\par
(ii) Similarly, given a smooth morphism of complex manifolds
$f:X_B\to B$, there is a natural relative $T_{X_B/B}$-deformation
parametrized by $P^m_B$, namely $\O_X\otimes_{\O_B}P^m_B$
(here $T_{X_B/B}$
denotes the Lie algebra of 'vertical' vector fields,
tangent to the fibres of $f$ .
We denote this deformation
by $P^m(X_B/B)$. Intuitively it represents
the family of $m-$th order deformations $f\inv(N_m(b)), b\in B.$
Since $T_{X_B/B}$ acts on
$\O_X$ by $\O_B-$linear derivations,
it follows that $P^m(X_B/B)$
is a relative deformation in the category of
$\O_B$-algebras.
\endremark
\ss

Now the construction of
universal deformations and related objects
extends in a
straightforward manner to the case of admissible $\g_B-$
deformations. Thus, there is a relative very symmetric product
$X<n>/B\overset f_n\to\longrightarrow B$
which is just the fibre product
$$X<n>\times_{B<n>}B<1>\to B<1>=B,$$
and on this we have a
relative Jacobi complex $J_m(\g_B/B)$ which has
a natural relative OS or comultiplicative structure, so that
$$\RR_m(\g_B/B):=
\O_B\oplus {\Cal Hom} (\R^0f_{m*}(J_m(\g_B/B)),\O_B)
=:\O_B\oplus\m_m(\g_B/B)$$
is a sheaf of $\O_B$-algebras of finite type as $\O_B-$module.
Moreover there is
a tautological morphic (comultiplicative) element
$$v_m\in\H^0(X<m>/B,J_m(\g_B)\otimes\m_m(\g_B/B)$$
and there is
correspondingly a tautological relative $\g_B-$ deformation
parametrized by
$\RR_m(\g_B/B)$, which we denote by $u_m/B.$ Under
suitable hypotheses, which we proceed to
state, $u_m/B$ and $v_m$ will be  universal.\par

Now the following result generalizes Theorem 3.1 above
and Theorem 0.1 of [Rcid], and can be proved similarly.
\proclaim{Theorem 9.2} Let $\g_B$ be an admissible \dgla over
$X/B$. Then
\item{(i)} to any isomorphism class of
relative $\g_B$-deformation parametrized by
an algebra $\SS$ of exponent $m$
there are canonically associated a morphic Kodaira-Spencer
element
$$\beta_m(\phi)\in\H^0(J_m(\g_B/B)\otimes\m_{\SS})$$
and a compatible
homomorphism of $\O_B-$algebras
$$\alpha_m(\phi):\RR_m(g_B/B)\to\SS;$$
conversely, any morphic element
$$\beta\in\H^0(J_m(\g_B/B)\otimes\m_{\SS})$$
induces a relative $\g_B$-deformation $\phi_m(\beta)$
parametrized by $\SS$;
\item{(ii)}if $\g_B$ has central sections then
there is an isomorphism of relative deformations
$$\phi\simeq\phi_m(\beta_m(\phi));$$
any two such isomorphisms differ by an element of
$$Aut(\phi)=H^0(\exp(\g_B^\phi\otimes\m_{\SS})).$$\endproclaim
\remark{Remarks 9.3}(i)\
As we have seen, there are nontrivial relative deformations
even if the fibres of $X_B\to B$ are points, in which case
$\RR_m(\g_B/B)=\O_B$ so $\alpha_m(\phi)$ certainly does not
determine $\phi$.\par\noindent
(ii)\ Note that in the above situation
$\RR(g_B/B)$ and $\SS$
are not necessarily $\O_B-$flat.\endremark
\remark{Example 9.4}
If $\SS$ is of exponent 1, i.e. $\m_{\SS}^2=0,$
then it is
easy to see directly that relative $\g_B$-deformations
parametrized by $\SS$ are in 1-1 correspondence with
$H^1(X,\g_B\otimes\m_S)$. The Kodaira-Spencer homomorphism
corresopnding to $\phi\in H^1(X,\g_B\otimes\m_S)$ is just
the corresponding map
$$(\R^1f_*(\g_B))^v\to\m_S.$$\endremark
 We might define a 'family of deformations
parametrized by $\SS$' to be a collection of isomorphism
classes of deformations over members of some open cover
of $B$, together with suitable gluing data over the overlaps;
this type of object is naturally classified by
$$H^0(\R^0f_{m*}(J_m(\g_B/B)^{vv}\otimes\m_{\SS})).$$
There is a natural map to this group from
$H^0(J_m(\g_B/B)\otimes m_{\SS})$, and assuming
$\R^0f_{m*}(J_m(\g_B/B)$ is locally free, this map
may be analyzed with the usual Leray spectral sequence,
which leads to the following result. First a definition.
We will say that a Lie algebra sheaf $\g_B$ as above has
{\it{relatively central sections}} if the image of
the natural map
$$f\inv f_*(\g_B)\to\g_B$$
is contained in the center of $\g_B$. Note that this condition is
stronger than saying that $\g_B$ has central sections, which
concerns the image of $H^0(X_B,\g_B)\to \g_B.$
\proclaim{Corollary 9.5} In the situation of Theorem 9.2,
assume additionally that $\g_B$ has relatively central sections,
that $\R^0f_{m*}(J_m(\g_B/B))$ is $\O_B$-locally free,
and that
$$H^i(f_*(\g_B)\otimes F)=0, \forall i>0,$$
for all coherent $\O_B-$modules $F$.
Then for any relative $\g_B$-deformation $\phi/\SS,$
$$\phi\simeq\alpha_m(\phi)^*(u_m)=u_m/B\otimes_{\RR_m(\g)}\SS.$$
In particular, relative $\g_B$-deformations are determined by
their associated Kodaira-Spencer homomorphisms.\endproclaim
\demo{proof} Our hypotheses imply that
$$H^i(\R^jf_{m*}J_m(\g_B/B))\otimes\m_{\SS})=0, \forall j<0,$$
so it suffices to apply the usual Leray spectral sequence to compute
$$H^0(J_m(\g_B/B)\otimes\m_{\SS})=
H^0(B,\R^0f_{m*}(J_m(\g_B/B))\otimes_{\SS}).\qed$$\enddemo
Note that the hypotheses of the Corollary are satisfied
provided first that $\R^0f_{m*}(J_m(\g_B/B)$ is locally free (i.e
$\g_B/B$ is 'relatively unobstructed'), and second, either
 $f_*(\g_B)=0$  or $B$ is an affine scheme ( provided all sheaves
 in question are coherent). In general however, a relative
 deformation cannot adequately by thought of as a family
 of isomorphism classes of deformations,
because gluing together isomorphism classes of deformation
is weaker than gluing together actual deformations.
\par
Finally we will show that the constructions
and results of \S 8 on connection
algebras carry over {\it{mutatis mutandis}} to the relative case.
Thus, suppose given an relative admissible pair $(\g_B,E_B)$ on
$X_B\overset f\to\to B$, with
soft $\O_B-$linear resolution $(\g^._B, E^._B)$, and assume given
a finite complex $F^.$ of free $\O_B$-modules of finite type
such that
$$H^j(F^.\otimes\C(b))
\simeq H^j(X_b, E_b),\ \forall j,\forall b\in B.$$
As is well known, such complexes $F^.$ always exist locally if
$f$ is a proper morphism of algebraic schemes and, as we shall see,
the final statement
will be essentially independent of the particular
complex $F^.$. Moreover, if $E_B$ is relatively $i-$ equicyclic
(i.e. $H^j(E_b)=0 \forall j\neq i$) we may assume
$F^j=0 \forall j\neq i, i-1.$ Then
there is a relative connection algebra
$$\k(\g_B,E_B):\ \g_B^.\to f\inv(^*F.)\otimes E_B^.$$
where $^*F^.=Hom^.(F^.,\O_B)$, which is still admissible
and acts on $E_B^.$, and the following relative
analogue of Theorem 8.2 holds.
\proclaim{Theorem 9.6} In the above situation,
assume additionally that
$\g_B$ has relatively central sections and that
$E_B$ is relatively
$i-$equicyclic. Then we have a class of isomorphisms
$$\R^if_*(M_m(\g_B,E_B))
\otimes_{\RR_m(\g_B/B)}\RR_m(\k(\g_B, E_B))\tag 9.3$$$$
\simeq \R^if_*(E_B)
\otimes_{\O_B}\RR_m(\k(\g_B, E_B))$$
any two of which differ by a map induced by an
element of $\Aut (u_m/B)$
where $u_m/B$ is
the universal relative deformation.\qed\endproclaim
\proclaim{Corollary 9.7}
In the situation of Theorem 9.6, assume moreover
that $f$ is a
smooth proper morphism of complex manifolds and that
for some $m\geq 2$ we have that
\item{(i)}
if $\phi_m$ is the relative deformation $P^m(E_B,\g_B)/B$
parametrized by $P^m_B$
(cf. Example 9.1(i)), then the associated
Kodaira-Spencer homomorphism
$$\alpha_m(\phi_m):\RR_m(\g_B/B)\to P^m_B$$
factors through $\RR_m(\k(\g_B, E_B))$;
\item{(ii)} $f\inv f_*(\g_B)$ acts on $E_B$ as scalars.\ss
Then the vector bundle $\R^if_*(E_B)$ admits a natural
projective connection.\endproclaim
\demo{proof} Set $G=\R^if_*(E_B)$. Then our assumptions give
an isomorphism of $P^m(G)$ and $G\otimes P^m_B$ globally defined
up to scalars.
For any $m\geq 2,$ this is equivalent to a projective connection.
\qed\enddemo
\vfill\eject
\heading 10. The Atiyah class of a deformation\endheading
Let $(\g_B,E_B)$ be an
admissible pair on $X_B/B$, $\SS$ a finite-length
$\O_B-$algebra, and $E^\phi$
an admissible $\g_B-$deformation parametrized
by $\SS.$ There is a corresponding deformation $\g^\phi$,
and clearly $\g^\phi$
is a Lie algebra acting on
$E^\phi$. We ignore momentarily the status of
$\E^\phi$ as a
deformation and just view it as a $\g^\phi-$module over
$X_\SS=X_B\times_B\Spec (\SS).$
Let $\Spec(\SS_1)$ be the
first infinitesimal neighborhood
of the diagonal in $\Spec (\SS)\times\Spec (\SS)$ with projections
$$p,q:\Spec(\SS_1)\to \Spec (\SS).$$
Then $p_*q^*E^\phi$ may be viewed as a first-order $\g^\phi$
deformation of $E^\phi$
and we let
$$AC(\phi)\in H^1(\g^\phi\otimes\m_{\SS_1})=
H^1(\g^\phi\otimes_\SS\Omega_{\SS/B})\tag 10.1$$
be the associated (first-order) Kodaira-Spencer class.\par
A cochain representative
for $AC(\phi)$ may be constructed as follows.
Let $$\phi\in\Gamma(\g^1)\otimes\m_\SS$$
be a Kodaira-Spencer
cochain corresponding to $E^\phi$ ,satisfying the
integrability condition
$$\p\phi=-\1/2 [\phi,\phi].$$
Let
$$d_\SS:\Gamma(\g^1)\otimes\m_\SS
\to \Gamma(\g^1)\otimes\Omega_{\SS/B}$$
be the map induced by exterior derivative on $\m_\SS$.
Set
$$\psi=d_\SS(\phi).\tag 10.2$$
Then
$$AC(\phi)=[\psi].$$
Note that
differentiating the integrability condition for $\phi$
yields
$$\p\psi=-[\phi,\psi].$$
Since $(\g^.,\p+\ad(\phi))$
is a resolution of $\g^\phi$, this means that
$\psi$ is a cocycle for $\g^\phi$.
\remark{Example 10.1} Let $E$ be a vector bundle
on $X_B$ with a $\g$-structure
as in Example 9.1.
Taking $\SS=P^1=\O_X\oplus\Omega_X$ as there, we get a
first-order relative $\g$-deformation $P^1(E,\g)$.
Note that in this case $\Omega_{\SS/B}=\Omega_X$ and its
$\SS-$module structure factors through $\O_X$. Thus
the Atiyah-Chern class
$$AC(P^1(E,\g))\in H^1(\g\otimes\Omega_X)$$
and it is easy to see that it coincides
with the usual
Atiyah-Chern class of the $\g$-structure $E$
which may be defined, e.g. differential-geometrically
in terms of a $\g$-connection
(and which reduces
to the usual Atiyah-Chern class if $\g=\gl(E)$, cf. [At]).
Indeed our good resolution in this case takes the form
$$E^0\otimes (\O_X\oplus\Omega_X)
\to E^1\otimes (\O_X\otimes\Omega_X)
\ldots$$
with differential
$$\left (\matrix \bp&\phi\\0&\bp\endmatrix\right ) $$
and note that in this case $\phi=\psi$ since
$\m_{\SS}=\Omega_{\SS}$. Assuming $E$ is endowed with a $\bp-$
connection, the parallel lift of a section $e$ of $E$
to $E^0\otimes (\O_X\oplus\Omega_X)$ is given by $(e,\nabla e)$
and consequently we have
$$\phi (e) = [\bp,\nabla](e).$$
Thus
$$\psi = [\bp,\nabla]\tag 10.1.1$$
In other words, for any section $v$ of $T_X$, holomorphic or not,
we have
$$\psi\neg v=[\bp,\nabla_v].\qed$$
\endremark
\remark {Example 10.2} Consider an ordinary first-order
deformation $\phi$ of a complex manifold $X$, corresponding to
an algebra $S$ of exponent 1. Suppose this deformation
comes from a geometric family
$$\pi:\Cal X\to Y$$
with $\Cal X, Y$ smooth, $S=\O_{Y,0}/\m_{Y,0}^2$.
Then it is easy to see that
$AC(\phi)$ corresponds to the extension
$$0\to T_X\to \D_\pi\to T_0Y\otimes\C_X\to 0$$
where $\D_\pi$ is the subsheaf of $T_{\Cal X}\otimes\O_X$
consisting of 'vector fields locally constant in the normal
direction', i.e. those derivations $\O_{\Cal X}\to\O_X$
that preserve the subsheaf $\pi\inv\O_Y\subset
\O_X.$\qed
\endremark
The last example suggests an interpretation of the
Atiyah class as an extension also in the general case.
To state this,
let $\phi$ be a relative deformation parametrized
by $\SS$ as above, and set
$$I=\Ann(\Omega_{\SS/B})\subset\SS, \SS'=\SS/I,
\phi'=\phi\otimes_{\SS}\SS'$$
and let $\Omega_{\SS/B}^{vv}$ denote the double dual as
$\SS'$-module. Note that
$$\Omega_{\SS/B}^{vv}=\Der_{\O_B}(\SS,\SS')^v$$
( dual as left $\SS'$-module). \par
We will also consider the analogous situation over
a formally smooth, complete noetherian augmented local
$\O_B$-algebra $\SS\w$ (which is
thus  locally a power series algebra over $\O_B$),
where of course dual means as (left) $\SS\w$-module.
\proclaim{Theorem 10.3} The image of $AC(\phi)$ in
$H^1(\g^{\phi'}\otimes\Omega_{\SS/B}^{vv})$ corresponds
to an extension of $\SS'$ modules
$$0\to\g^{\phi'}\to\D(\phi)\to
f\inv \Der_{\O_B}(\SS,\SS')\to 0\tag 10.3$$
and there is a natural action pairing
$$\D(\phi)\times E^\phi\to E^{\phi'}.$$
Moreover, if $\phi\w$ is a formal deformation parametrized
by a formally smooth
$\O_B$-algebra $\SS\w$, then the image of
$AC(\phi\w)$ in
$H^1(\g^{\phi\w}\otimes\Omega_{\SS\w/B}^{vv})$
corresponds to an extention of $\SS\w$-Lie algebras
$$0\to\g^{\phi\w}\to\D(\phi\w)\overset\nu\to\to
f\inv T_{\SS\w}\to 0\tag 10.4$$
where $T_{\SS\w}=\Der_{\O_B}(\SS\w,\SS\w)$ and
$\D(\phi\w)$ acts on $E^{\phi\w}$ satisfying
the rule
$$d(f.v)=f.d(v)+\nu(d)(f).v,
\forall d\in\D(\phi\w),
f\in\SS\w,v\in E^{\phi\w}\tag 10.5$$
\endproclaim
\demo{proof} For brevity we shall work out the formal case,
the artinian case being similar. As usual we let $(\g^.,E^.)$
be a soft (dgla,dg module) resolution of $(\g,E)$; also let
$(C^.,\p)$ be a soft resolution of $f\inv\O_B$, and
note that $\g^.$ is a $C^.$-module. Then clearly
$\D(\phi\w)$, i.e.
the extension corresponding to $AC(\phi\w)$ is resolved
by the complex
$$\D^.(\phi\w)=\g^.\otimes\SS\w\oplus C^.\otimes T_{\SS\w} $$
with differential given by the matrix
$$\left (\matrix \p+\phi\w&\psi\w\\
0&\p \endmatrix\right ) $$
where $\psi\w=d_{\SS\w}(\phi\w)$ as in (10.2), which defines
in an obvious way a map
$C^i\otimes T_{\SS\w}\to\g^{i+1}\otimes \SS\w$.
\par
Now we claim that $\D^.(\phi\w)$ is a dgla: indeed since
$\g^.\otimes\SS\w$ and $ T_{\SS\w}\otimes C^.$ with the induced
differentials are
clearly dgla's (in the latter case, the bracket is
induced by that of $T_{\SS\w}$), and
$ T_{\SS\w}\otimes C^.$ acts on $\g^.\otimes\SS\w$
via the action of
$ T_{\SS\w}$ on $\SS\w$ and the $C^.-$module
structure of $\g^.$,
it suffices to show that
$\psi\w$ is a derivation, which is essentially obvious:
$$\phi\w([v_1,v_2])=
[v_1,v_2](\phi\w)=v_1(v_2(\phi\w))-v_2(v_1(\phi\w))$$
$$v_1(\psi\w(v_2))-v_2(\psi\w(v_1)).$$
Now since
$\D^.(\phi\w)$ is a dgla, the fact that it acts on $E^{\phi\w}$
essentially follows from the fact
that the differential of $\D^.(\phi\w)$
is just
commutator with the differential on the resolution of
$E^{\phi\w}$,
i.e. $\p+\phi\w$. To check the latter, it is firstly clear
on the $\g^.\otimes\SS\w$ summand; for the other summand,
take $v\in T_{\SS\w}\otimes C^.$. Then
$$[v,\p+\phi\w]=[v,\p]+[v,\phi\w]=\p(v)+\psi\w(v).$$
This shows that the obvious term-by-term pairing induces
a pairing of complexes
$$\D^.(\phi\w)\times (E^.,\p+\phi\w)\to (E^.,\p+\phi\w),$$
 whence a pairing
$\D(\phi\w)\times E^{\phi\w}\to E^{\phi\w}$;
that this is in fact
a Lie action is clear from the fact that the corresponding
assertion holds term-by-term. This completes the proof.
\qed
\enddemo
\vfill\eject
\heading 11. Vector bundles on manifolds:
the action of base motions\endheading

In this section we
go back to the
situation considered in \S 6, with a locally fine moduli
space $\M$
with associated Lie algebra $\tilde{\g}$ on $X\times \M$.
We assume additionally
that $X$ is a compact complex manifold  and
$\tilde{\g}$ is an $\O_X-$Lie algebra
$\g$ acting $\O_X$-linearly.  We assume that
$$R^0p_{\M*}(\tilde{\g})=0.\tag 11.1$$
For convenience, we shall also assume that
$$R^2p_{\M*}(\tilde{\g})=0,\tag 11.2$$
which in particular implies that $\tilde{\g}$ is (relatively)
unobstructed,
so that $\M$ is smooth (it seems reasonable that similar
results can be
obtained assuming only the unobstructedness).
Of course, condition (11.2) holds automatically
when $X$ is a Riemann surface.\par
Since $\M$ is in a sense a functor of $X$,
it seems intuitively
plausible that a motion- say an infinitesimal
 motion, i.e. global holomorphic vector field on $X$-
  should induce a similar motion of $\M$.
In this naive form this intuition seems of little use per se,
since in cases
of interest $X$ will not admit any {\it global}
 holomorphic vector fields while {\it local}
vector fields have no obvious relation to $\M$.
But there is another, more 'global' way to represent
the Lie algebra $T_X$ of holomorphic vector fields on $X$,
namely via the Dolbeault algebra $A^.(T_X).$
Then the 'induced motion' idea
suggests that there should be (something like) a map
$$\Sigma: A^.(T_X)\to A^.(T_\M).\tag 11.3$$
Since $\Sigma$, at least in some sense, sends a motion
of $X$ to the induced motion
of $\M$,it should be a dgla homomorphism.
Now, at least on cohomology,
a map as in (11.3) exists: it is none other
than 'cap product with the
Atiyah class of the universal bundle'
which indeed is
given essentially just by differentiating a cocycle
defining this
universal bundle with respect to the given vector field,
then pushing down to $\M$.
The upshot, then, is that a suitable
version of the map $\Sigma$ ought to be a Lie homomorphism, i.e.
compatible with brackets
(as well as, of course, the differential).
This is what
we aim to show in this section. As one might expect,
this fact is
important in relating deformations of $X$ and $\M$.\par
The map $\Sigma$ is defined as follows. Let
$$\psi\in\Gamma (\gt^1)\otimes\Omega_{X\times\M}$$
be a representative of the Atiyah class
$$[\psi]=AC(P^1(\gt,\gt))\in
H^1(\gt\otimes\Omega_{X\times\M}$$
of the 1st jet of $\gt$
over $X\times\M$ (cf. Example 10.1).
We replace $T_X$ by its Dolbeault resolution $A^.(T_X)$
(where $\cdot =(0,\cdot)$), truncated beyond
degree 2
(which doesn't affect the deformation theory),
 and define a map
$$\Sigma_0:A^.(T_X)\to A^{.+1}_{X\times\M}(\gt),$$
$$\Sigma_0 (v) = \psi\neg v,\ \ v\in A^.(T_X)\tag 11.4$$
where $\neg$ denotes interior multiplication or contraction.
 Since $\psi$ is $\bp-$closed,
clearly $\Sigma_0$ commutes with $\bp$. On the other hand
our assumptions (11.1), (11.2) plus
the fact that $\M$ is locally a fine moduli imply that the
analogous map
$$\Sigma_1:A^.(T_\M)\to A^{.+1}_{X\times\M}(\g),$$
$$\Sigma_1(v)=\psi\neg v,\ \ v\in A^.(T_\M)\tag 11.5$$
is a quasi-isomorphism,
so we get a map in the derived category
$$\Sigma =\Sigma_1\inv\circ\Sigma_0\tag 11.6$$
Our main result concerning $\Sigma$ is the following
\proclaim{Theorem 11.1}
$\Sigma$ is a dgla homomorphism, i.e. is compatible
with brackets.\endproclaim
\demo{proof}
It clearly suffices to prove that if $v_1,v_2\in A^0(T_X)$,
$$v_i=\sum a_{i,j}\p/\p z_j$$
are two type-(1,0) vector fields (not necessarily holomorphic),
then
$$[\Sigma(v_1),\Sigma(v_2)]=\Sigma([v_1,v_2])\tag 11.7$$
To show that the two sides of (11.7) agree it suffices to
check they agree pointwise at each point of $\M$. To
this end we will use the recipe of \S 6 to compute the LHS.\par
So let us fix a point $z$ of $\M$, corresponding to a particular
pair $(\g,E)$, and fix a
$ g-$connection of type $\bp$ on $E$ and $\g$. Then
first of all, it is clear by Example 10.1 that the
'value' of $\Sigma (v)$ at any point $w\in\M$ is given by
$$\Sigma (v)|_w= [\nabla_v,\bp_w]\tag 11.8$$
where $\bp_w$ is the $\bp$ operator corresponding to $w$.
Next, consider the restriction of $\Sigma (v)$ on the first
infinitesimal
neighborhood $N_1(z)$ of $z$ in $\M$. By universality,
we may identify
the restriction of $\g$ on $X\times N_1(z)$ with
$\g^\phi$,
the first-order infinitesimal $\g$-deformation of $\g$,
and likewise for $E$.
Let $(\phi_i\in\Gamma (\g^1))$ be a lift of a basis
of $H^1(\g)$,
and $(\phi_i^*)$ a lift of a dual basis. Now the
prolongation of $\bp_z$
in the direction corresponding to $\phi_i$
is obviously given by $\bp_z+\phi_i$, hence may write
$$\bp|_{N_1(z)}=\bp_z+\sum\phi_i^*\otimes\phi_i.$$
Therefore by (11.8) we have
$$\Sigma (v)|_{N_1(z)}=
[\bp_z,\nabla_v] + \sum\phi_i^* [\phi_i,\nabla_v]$$
Note that $\sum\phi_i^* [\phi_i,\nabla_v]$ is just the cobracket
$^t\br (\nabla_v).$
Now by elaboration 6.2 we compute:
$$[\Sigma (v_1),\Sigma (v_2)]|_z=
<^t\br (\nabla _{v_1}),[\bp_z,\nabla_{v_2}]>-
<^t\br (\nabla _{v_2}),[\bp_z,\nabla_{v_1}]>$$
$$=[[\bp_z,\nabla_{v_2}],\nabla_{v_1}]-
[[\bp_z,\nabla_{v_1}],\nabla_{v_2}].$$
Applying the Jacobi identity to the first term yields
$$[\Sigma (v_1),\Sigma (v_2)]|_z=
-[[\nabla_{v_2},\bp_z],\nabla_{v_1}]-
[[\nabla_{v_1},\nabla_{v_2}],
\bp] -[[\bp_z,\nabla_{v_1}],\nabla_{v_2}]=
[\bp_z,[\nabla_{v_1},\nabla_{v_2}]].$$
But as our connection of of $\bp$ type, its curvature is
of type $(1,1)$, while $v_1,v_2$ are of type $(1,0)$, hence
$$[\nabla_{v_1},\nabla_{v_2}]=\nabla_{[v_1,v_2]}.$$
Consequently we have
$$[\Sigma (v_1),\Sigma (v_2)]|_z=
[\bp_z,\nabla_{[v_1,v_2]}]=\Sigma ([v_1,v_2])|_z.$$
Therefore (11.7) holds and the proof is complete.\qed

\enddemo
\vfill\eject
\heading 12. Vector bundles on Riemann surfaces: refined action by
base motions and Hitchin's connection\endheading
Our purpose here is to refine the results of the previous section,
in the case
where $X$ is 1-dimensional, by constructing a lift of $\Sigma$
to another dgla associated to $\M$. We continue with the notations of
that section; in particular, $\M$ is a locally fine moduli space
associated to a dgla sheaf $\gt$ on $X\times\M$, and we also fix a
$\gt-$ deformation $\Et$ on $X\times\M$, such that
$$R^ip_{\M*}(\Et)=0,\ \ i\neq 1,$$
and consequently
$$G:=\bigwedge\limits^{top}R^1p_{\M*}(\Et)$$
is an invertible sheaf on $\M$.\par
 We note that $G$ itself may be
realized as the (sole nonvanishing)
direct image of a suitable $\gt-$deformation,
as follows. Note that
$$\gt_r:=\pi_{r*}p_1^*(\gt),$$
where
$\pi_r:X^r\times\M\to X<r>\times\M ,\ \ X^r\times\M\to X\times\M$
are natural
projections, naturally has the structure of dgla sheaf acting
on $\lambda^rE$, and clearly
$$G=R^rp_{\M*}(\lambda^rE)$$
with all other derived images being zero.
There is a pullback map
$$Rp_{\M*}(\gt)\to Rp_{\M*}(\gt_r)$$
which is compatible with brackets and induces
isomorphisms on $R^0$ and $R^1$. Choosing a fixed base-set
$\{ x_1,...,x_{r-1}\}\in X<r-1>$ yields an embedding
$X\to X<r>$ which induces a splitting of the pullback map,
showing that this map is injective on $R^2$. It follows that
we have a natural isomorphism
$$R^0p_{\M*}(\gt)\to R^0p_{\M*}(\gt_r).$$
Hence we may view $G$ as the direct image of a $\gt$-deformation.
\par
As in Example 1.1.2 C, \S 1, we may consider the heat atom
$\D^{2/1}(G)$ associated to the $\O_\M-$module
$G$, which is the pair
$$\D^1(G)\to \D^2(G).$$
Note that since $G$ has rank 1, $\D^{2/1}(G)$ is equivalent as
complex on $\M$ to $\Sym^2T_\M[-1]$, which is thus endowed with
a Lie bracket.
Also, $\D^{2/1}(G)$ is obviously equivalent to the pair
(Lie atom) $\D^1(G)/\O\to \D^2(G)/\O.$
 As we have seen, $\D^i(G)$ may be naturally identified
with the direct image
of $J_i(\gt, F^*)\otimes F$ where $F=\lambda^r(E)$
hence $\D^{2/1}(G)$ ,
i.e. $\Sym^2T_\M[-1]$ is the direct image of
$\lambda^2(\gt)[1]$.
\par
We now assume $X$ is of dimension 1 and that $\M$ is the
global fine moduli
space ${\Cal SU}^r_X(d)$ or  ${\Cal SU}^r_X(L)$
of vector bundles of rank $r$
and fixed determinant $L$ on $X$, where $L$ is
a line bundle of degree $d$.
We assume temporarily that $(d,r)=1$ (the
modifications needed to handle the general case
will be indicated later. As is well known [NaRam],
the assumption $(d,r)=1$ implies that $\M$ is a
fine moduli, in particular a locally fine
moduli space associated to the Lie algebra
$\gt=\frak{sl} (E)$, in the sense of \S 6.
By Proposition 4.8, the map $\Sigma$ in this case factors
through $\lambda^2(\gt)[1]$,
it follows that $\Sigma$ factors through a map
$$\Omega: T_X\to \D^{2/1}(G).$$
\proclaim{Theorem 12.1}
$\Omega$ is a Lie homomorphism.\endproclaim
\demo{proof}
Recall that we are identifying $T_X$ with the dgla $A^.(T_X),$
which exists in
degrees 0,1. In degree 0, $\D^{2/1}(G)$ can be identified
with $\D^1(G)/\O\simeq T_\M$, so the homomorphism
property is just Theorem 10.1.
Therefore it just
remains  to prove the homomorphism property in degree 1.
For any $v\in T_X,$ write
$$\Omega (v)=(A(v), B(v)),$$
with $A(v)\in \D^1(G)/\O, B(v)\in \D^2(G)/\O.$
Then what has to be shown is that
for any $v_0\in A^0(T_X), v_1\in A^1(T_X),$ we have
$$ B([v_0,v_1])=<A(v_0),B(v_1)-<A(v_1),B(v_0)>.\tag 12.1$$
Now firstly, $B(v_0)=0$ since $B$ lowers degree by 1.  Next, since
$[v_0,v_1]$ is automatically $\bp-$closed, we have
$$\bp B([v_0,v_1])= A([v_0,v_1]) = [A(v_0),A(v_1)]$$
the last equality by Theorem 11.1.
Again because $v_1$ is $\bp-$closed, we
have
$$A(v_1)=\bp B(v_1).$$
The upshot is that both sides of (12.1) have the same $\bp$,
hence their difference
yields a global holomorphic section of $\D^2(G)/\O$ over $\M$.
However, it is well known that
$\D^2(G)/\O$ has no nonzero sections:
indeed this follows from Hitchin's result
that the coboundary map
$$H^0(\Sym^2T_\M)\to H^1(T_\M)$$
is injective, plus the fact that $H^0(T_\M)=0$ (cf. [NaRam]).
This completes the proof.
\qed\enddemo
Now consider the diagram
$$\matrix \D^1(G)/\O&\to&\D^2(G)/\O\\
\downarrow&&\downarrow\\
G\otimes^*G/\O I&=&G\otimes ^*G/\O I\endmatrix$$
where the vertical arrows are induced by the action of $\D^1(G)$
and $\D^2(G)$ on $^*G$ and $I$ is the identity in $G\otimes ^*G$.
This diagram itself may be considered
a dgla quasi-isomorphic to $\D^{2/1}(G).$
And of course the left column is quasi-isomorphic to
$\k (\D^1(G),G)$ (cf. \S 8).
Consequently, we have a Lie homomorphism
$$\D^{2/1}(G)\to\k (\D^1(G),G).$$
Composing this with $\Omega$ above, we get a Lie homomorphism
$$\omega: T_X\to \k (\D^1(G),G).$$
It follows easily from this that over the deformation space of pairs
$(X,L)$ there is a canonical
local
trivialization or connection on
 the projective bundle associated to $H^0(G)$,
which is the main result of Hitchin [Hit] (see also
[BryM], [F],  [Ram], [Sun],
[TsUY], [vGdJ], [WADP] and references therein;
the connection is sometimes called the Hitchin or
Knizhnik-Zamolodchikov
connection):

\proclaim{Corollary 12.2} Let $Y$ be any manifold parametrizing pairs
$(X,L)$ where $X$ is a compact Riemann surface of genus $g\geq 3$
and $L$ is a line bundle of degree $d$ on $X$, and let
$\SH$ be the vector bundle on $Y$ with fibre
$H^0({\Cal SU}^r(X,L), G)$. Assume $d,r$ are relatively prime.
Then there is a canonical projective connection on $\SH$.
\endproclaim
\demo{proof} We have a family of smooth curves $X_Y/Y$
and a family of associated moduli spaces which we denote by
$\M_Y/Y$, and there is a commutative diagram of $\O_Y$-algebras
and homomorphisms:
$$\matrix R_m(T_{X/Y}/Y)&\to&P^m_Y\\
\uparrow&\nearrow&\\
R_m(T_{\M_Y}/Y)&&\endmatrix$$
where the vertical homomorphism is induced by $\Sigma.$
This diagram represents the intuitive fact that we have
a family of $m-$th order deformations of fibres $X_y$
and $\M_y$ for $y\in Y$ (cf. Example 9.1(ii)).
As we have seen in Theorem 12.1,
the map induced by $\Sigma$ factors through
$$\tilde{R}=R_m(\k(\D^1(G), G)).$$ The
module $P^m(\SH)$ comes by extension of scalars
from an analogous module over
$\tilde{R}$ which by Corollary
7.2 is isomorphic (up to scalars) to $\SH\otimes_{\O_B}
\tilde{R}$.
Hence as $\O_Y$-modules, $P^m(\SH)$ and $P^m_Y\otimes\SH$
are isomorphic up to scalars, so there is
a projective connection (cf. Corollary 9.7).
\qed
\enddemo
Now we will indicate the extension of this result to the case
where $d$ and $r$ have a common factor, so that we have only
a coarse moduli space without a universal family. Fixing $r,d$,
let $U^s\subset SU^r(X,d)$ be the subset corresponding to stable
bundles. As is well known, under our assumptions $SU^r(X,d)$ is
normal and projective and the complement of $U^s$ has codimension
$>1$, hence for any line bundle $F$ on $SU^r(X,d)$  the restriction
map $$H^0(F,SU^r(X,d) )\to H^0(F,U^s)$$ is an isomorphism. Now by
construction (see [NaRam], [Sesh], [VLP]),
there is a finite collection $\UU$ of locally fine
moduli spaces $U_\alpha$, with corresponding rank-$r$ universal
bundles $\Et_\alpha$ on $X\times U_\alpha$, such that the images
of the natural maps.
$$f_\alpha:U_\alpha\to U^s$$
form a covering.
 We may further assume that each $U_\alpha$ is
affine and Galois over its image in $U^s$, and that the
collection $(U_\alpha,f_\alpha)$ is 'Galois-stable' in the sense
that for each deck transformation $\rho, (U_\alpha,f_\alpha\circ
\rho)$ is also in the collection.
Now set
$$U_{\ab}:=U_\alpha\times_{U^s} U_\beta$$
and likewise for triple products etc. Let
$$p_\alpha:=1_X\times f_\alpha:X\times U_\alpha\to X\times U^s,$$
Also let
$$p_{\ab,\alpha}:X\times U_{\ab}\to X\times U_\alpha$$
be the obvious projection, and let
$$p_{\ab}:X\times U_{\ab}\overset {p_{\ab,\alpha}}\to
\to X\times U_\alpha\overset{p_\alpha}\to\to
X\times U^s$$ be the composite, and again likewise for higher products.
Note that
for any coherent sheaf $F$ on $X\times
U^s$, we may form a \v{C}ech-type complex
(of sheaves)
$$\check{C}(\UU, F):\bigoplus\limits_\alpha p_\alpha^*F
\to\bigoplus\limits_{\alpha,\beta}p_{\ab}^*F\to\cdots$$
and our saturation condition ensures that
the cohomology of $F$-$ H^0$ included-
may be computed from the hypercohomology of this complex,
in other words $\check{C}(\UU, F)$ is quasi-isomorphic to $F$, i.e.
to its \v{C}ech complex with respect to an ordinary
cover of $X\times U^s$ (thus ' \'etale
cohomology coincides with ordinary cohomology for coherent
sheaves'). Of course in our case the problem is
that we don't have an actual universal bundle $E$,
defined as a sheaf over all of $X\times U^s$ (this is
a result of Nori, cf. [Sesh]). However,
we shall see that we can still define a complex
to play the role of $\check{C}(\UU, E)$
for a universal bundle $E$, and the foregoing discussion
shows that this is 'good enough' at least for cohomology.
\par Note next that up to
shrinking our cover, we may assume we have isomorphisms
$$\sigma_{\beta\alpha}: p_{\ab,\alpha}^*\Et_\alpha
\to p_{\ab,\beta}^*\Et_\beta.$$

Indeed the sheaf
$$p_{U_{\alpha\beta }*}(\Cal Hom (
p_{\ab,\alpha}^*\Et_\alpha, p_{\ab,\beta}^*\Et_\beta))$$
is invertible by stability, hence after shrinking may
be assumed trivial, and a nonvanishing section of it yields the required
isomorphism. There is obviously no loss
of generality in assuming that
$$\sigma_{\beta\alpha}=\sigma_{\ab}\inv .$$
 Now note that over a triple product
$U_{\alpha\beta\gamma}:=U_\alpha\times U_\beta\times U_\gamma$,
the map
$$\sigma_{\gamma\alpha}
\inv\circ\sigma_{\gamma\beta}\circ\sigma_{\beta\alpha}
\in\Aut (\Et_\alpha|_{U_{\alpha\beta\gamma}})$$
must, for the same reason, be a scalar. Consequently,
$\sigma_{\gamma\alpha}
\inv\circ\sigma_{\gamma\beta}\circ\sigma_{\beta\alpha}$ induces
the $identity$ on
$$\gt_\alpha|_{U_{\alpha\beta\gamma}}
 =\sll (\Et_\alpha)|_{U_{\alpha\beta\gamma}}
 \subset  \Et_\alpha\otimes\Et_\alpha^*|_{U_{\alpha\beta\gamma}}.$$
 Consequently,
we may form a complex which may be considered the '\v{C}ech complex'
for $\gt$ with respect to the \'etale covering $\UU:=(U_\alpha)$:
namely the complex with sheaves
$$\bigoplus\check{C}(\gt, \UU)_{\alpha\beta\gamma...}
=\bigoplus\gt_\alpha|_{U_{\alpha\beta\gamma...}}
:=\bigoplus p_{\alpha\beta\gamma...,\alpha}^*\gt_\alpha,$$
each of which we identify  with its own Dolbeault or \v{C}ech
complex (using some affine covering of $X$),
and whose differentials are constructed as usual from the
pullback maps
$$r_{\alpha\beta\gamma...,\alpha\beta\gamma...\epsilon...}:
\gt_\alpha|_{U_{\alpha\beta\gamma...}}\to
\gt_\alpha|_{U_{\alpha\beta\gamma...\epsilon...}}$$
and from maps
$$r_{\alpha\beta\gamma...,\epsilon\alpha\beta\gamma...}:
\gt_\alpha|_{U_{\alpha\beta\gamma...}}\to
\gt_\epsilon|_{U_{\epsilon\alpha\beta\gamma...}}$$
given by restriction to
$\gt_\alpha|_{U_{\alpha\beta\gamma...\epsilon...}}$
followed by the isomorphism
$$\gt_\alpha|_{U_{\alpha\beta\gamma...\epsilon...}}\to
\gt_\epsilon|_{U_{\epsilon\alpha\beta\gamma...}}$$
induced by $\sigma_{\epsilon\alpha}$. By the above, these
indeed form a complex, and this
complex automatically inherits the structure of
a dgla from $\gt$. For the purposes of our constructions,
this complex may be taken as a substitute for $\gt$ itself.
Moreover, since the adjoint action of $\gt$ on itself is
faithful, we may take $\gt$ as a substitute for the universal
$\gt$ deformation $\Et$. \par
Now of course the theta-bundle $\theta$ itself and its powers
such as $F= \det H^1(\gt)$  of course exist
as actual line bundles on $U^s$, and all the auxiliary complexes
we need are derived from $\gt$ and $F$. Note that for any
line bundle $L$ we have a natural
 isomorphism of Lie algebras
$$\D^1(L)\simto \D^1(L^k),$$
given by the formula
$$D(s_1\cdots s_k)=\sum s_1\cdots D(s_i)\cdots s_k,$$
where $D$ and $s_1,...,s_k$ are local sections of
$D^1(L)$ $ L$ respectively.
Consequently, we may identify $\D^1(\theta^k)$ and $D^1(F)$
as Lie algebras. Hence all of our constructions go  through
in this context and establish the flatness of the connection.
\proclaim{Corollary 12.2 bis} The conclusion of Corollary 12.2
holds for all $d,r$.\qed\endproclaim
\vfill\eject
\heading References\endheading
\item{[At]}
Atiyah, M.F.: Complex analytic connections in fibre bundles.
Transac. Amer. Math. Soc. 85 (1957), 181-207.
\item{[BeK]}
Beilinson, A., Kazhdan, D.: Projective connections.\nuline
\item{[BryM]} Brylinski, J.-L., McLaughlin, D.:
Holomorphic quantization
and unitary representations of the Teichm\"uller
group. In: Lie theory and geometry, ed. Brylinski et al., Prog.
Math. 123 (1994).
\item{[DrNa]} Drezet, J.-M., Narasimhan, M.S.: Groupe de Picard
des vari\'et\'es de fibr\'es semistables sur les courbes
alg\'ebriques. Invent. math. 97, 53-94 (1989).
\item{[F]} Faltings, G.:
Stable G-bundles and projective connections. J. Alg. Geom.
2, 507-568 (1993).
\item{[Hit]} Hitchin, N. J.: Flat connections and geometric
quantisation. Commun.
math. Phys. 131, 347-380 (1990).
\item{[NaRam]} Narasimhan, M.S., Ramanan, S.: Deformations
of the moduli space of vector bundles over an algebraic curve.
Ann. Math. 101 (1975), 39-47.
\item{[Ram]}
Ramadas, T.R.: Faltings' construction of the K-Z connection.
Commun. Math. Physics 196, 133-143 (1998).
\item{[Rcid]}
Ran, Z.: Canonical infinitesimal deformations. J. Alg.
Geometry 9, 43-69 (2000).
\item{[Ruvhs]}Ran, Z.: Universal variations of Hodge structure.
Invent. math. 138, 425-449 (1999).
\item{[Sesh]}Seshadri, C.S.: Fibr\'es vectoriels sur les courbes
alg\'ebriques. Ast\'erisque 96 (1982).
\item{[Sun]} Sun, X.: Logarithmic heat projective operators.
Preprint.
\item{[TsUY]} Tsuchiya, A., Ueno, K., Yamada, Y.: Conformal
field theory on universal family of stable curves with
gauge symmetries. Adv. Stud. Pure math. 19,459-566 (1989).
\item{[vGdJ]} van Geemen, B., de Jong, A.J.:
On Hitchin's connection. J. Amer. Math. Soc. 11, 189-228 (1998).
\item{[VLP]} Verdier, J.-L., Le Potier, J., edd.: Module
des fibr\'es stables sur les courbes alg\'ebriques. Progr. in
Math. 54 (1985), Birkh\"auser.
\item{[WADP]} Witten, E., Axelrod, S., Della Pietra, S.:
Geometric quatization of Chern-Simons gauge theory. J. Diff.
Geom. 33, 787-902 (1991).
\item{[Welters]}
Welters, G.: Polarized abelian varieties and the heat equations.
Comp, math. 49, 173-194 (1983).

\enddocument